\documentclass[11pt, reqno, psamsfonts]{amsart}
\pdfoutput=1

\usepackage{amssymb}
\usepackage{amsthm}
\usepackage{amsmath}
\usepackage{latexsym}
\usepackage[T1]{fontenc}
\usepackage[utf8]{inputenc}
\usepackage[russian, french, english]{babel}

\usepackage{graphicx}
\usepackage{wrapfig}
\usepackage[justification=centering, labelfont=bf]{caption}
\usepackage{mathtools}
\usepackage[hidelinks]{hyperref}
\usepackage{tikz}
\usepackage{amsbsy}
\usepackage{enumitem}
\usepackage{mathrsfs}
\usetikzlibrary{shapes,snakes}
\usetikzlibrary{arrows.meta}
\usetikzlibrary{decorations.pathreplacing}
\usepackage{float}
\usepackage{multicol}
\usepackage{bbm}
\usepackage{bm}
\usepackage{mathabx}

\usepackage{lmodern}
\usepackage{titlesec}

\usepackage[spacing=true,kerning=true,babel=true,tracking=true]{microtype}

\usepackage[shortcuts]{extdash}

\usepackage[foot]{amsaddr} 

\usepackage{framed}

\usepackage[left=1in,right=1in,top=1in,bottom=1in,bindingoffset=0cm]{geometry}

\usepackage[backend=biber, style=alphabetic, sorting=nyt, maxnames=100]{biblatex}
\title[Ergodic Theorems and Pointwise Versions of The Ab\'ert--Weiss Theorem]{Ergodic Theorems for the Shift Action and Pointwise Versions of The Ab\'ert--Weiss Theorem}
\date{}
\author{Anton~Bernshteyn}
\address{Department of Mathematics, University of Illinois at Urbana--Champaign, IL, USA and Department of Mathematical Sciences, Carnegie Mellon University, Pittsburgh, PA, USA}
\email{bernsht2@illinois.edu; abernsht@math.cmu.edu}
\thanks{This research is supported in part by the Waldemar J., Barbara G., and Juliette Alexandra Trjitzinsky Fellowship.}

\newtheoremstyle{bfnote}%
{}{}%
{\slshape}{}%
{\bfseries}{\bfseries.}%
{ }%
{\thmname{#1}\thmnumber{ #2}\thmnote{ \ep{\normalfont{}#3}}}

\newtheoremstyle{defbfnote}%
{}{}%
{}{}%
{\bfseries}{.}%
{ }%
{\thmname{#1}\thmnumber{ #2}\thmnote{ (#3)}}

\newtheoremstyle{claim}%
{}{}%
{\slshape}{}%
{\itshape}{.}%
{ }%
{\thmname{#1}\thmnumber{ #2}\thmnote{ \ep{\normalfont{}#3}}}

\theoremstyle{bfnote}

\newtheorem{theo}[equation]{Theorem}
\newtheorem{prop}[equation]{Proposition}
\newtheorem{lemma}[equation]{Lemma}
\newtheorem{corl}[equation]{Corollary}

\newtheorem*{claim*}{Claim}
\newtheorem*{corl*}{Corollary}

\theoremstyle{claim}

\newcounter{ForClaims}[section]
\newtheorem{claim}{Claim}[ForClaims]

\newcommand*{\myproofname}{Proof}
\newenvironment{claimproof}[1][\myproofname]{\begin{proof}[#1]\renewcommand*{\qedsymbol}{\(\dashv\)}}{\end{proof}}

\theoremstyle{definition}
\newtheorem{defn}[equation]{Definition}
\newtheorem*{defn*}{Definition}

\newtheorem{remk}[equation]{Remark}

\newtheorem{remks}[equation]{Remarks}
\newtheorem*{exmp*}{Example}

\newtheorem{prob}[equation]{Problem}

\theoremstyle{remark}
\newtheorem*{ques*}{Question}
\newtheorem*{remk*}{Remark}

\renewcommand{\qedsymbol}{$\blacksquare$}

\newcommand{\0}{\varnothing}
\newcommand{\set}[1]{\{#1\}}

\newcommand{\dom}{\mathrm{dom}}

\newcommand{\supp}{\mathrm{supp}}
\newcommand{\Free}{\operatorname{Free}}
\newcommand{\acts}{\curvearrowright}
\newcommand{\N}{{\mathbb{N}}}
\newcommand{\Z}{\mathbb{Z}}
\newcommand{\F}{\mathbb{F}}
\newcommand{\R}{\mathbb{R}}
\renewcommand{\C}{\mathbb{C}}
\newcommand{\Q}{\mathbb{Q}}
\newcommand{\id}{\operatorname{id}}
\renewcommand{\epsilon}{\varepsilon}
\renewcommand{\phi}{\varphi}
\renewcommand{\theta}{\vartheta}
\newcommand{\pr}{\mathsf{p}}
\newcommand{\uni}{\upsilon}

\newcommand{\code}[1]{\pi_{#1}}

\newcommand{\E}{\mathbb{E}}
\newcommand{\M}{\mathbb{M}}

\renewcommand{\leq}{\leqslant}
\renewcommand{\geq}{\geqslant}
\renewcommand{\preceq}{\preccurlyeq}

\newcommand{\symdif}{\bigtriangleup}

\newcommand{\fins}[1]{[#1]^{<\infty}}
\newcommand{\finf}[2]{[#1 \to #2]^{<\infty}}
\renewcommand{\G}{\Gamma}
\renewcommand{\Prob}{\operatorname{Prob}}
\newcommand{\defeq}{\coloneqq}
\newcommand{\rest}[2]{{{#1}\vert_{#2}}}
\newcommand{\emphd}[1]{{\fontseries{b}\selectfont\textsf{#1}}}
\newcommand{\B}{{\mathscr{B}}}

\newcommand{\Def}{\mathrm{Def}}
\newcommand{\Stab}{\mathrm{St}}
\renewcommand{\P}{\mathbb{P}}
\newcommand{\p}{\pi}
\newcommand{\Ind}{\mathrm{Ind}}

\newcommand{\dist}{{\mathfrak{d}}}
\newcommand{\Diff}{\mathrm{d}}
\newcommand{\Nbhd}{\mathrm{N}}

\newcommand{\LG}{{\bm{\lambda}}}

\newcommand{\bemph}[1]{{\normalfont#1}} 
\newcommand{\ep}[1]{\bemph{(}#1\bemph{)}} 

\newenvironment{coolproof}[1][Proof]{\begin{proof}[\textsc{#1}]}{\end{proof}}

\usepackage{xspace}
\newcommand{\pmp}{$\text{p.m.p.}$\xspace}

\numberwithin{equation}{section}

\renewcommand{\thesubsection}{\arabic{section}.\Alph{subsection}}

\usepackage{etoolbox}

\titleformat{\section}[block]{\scshape\filcenter}{\thesection.}{1ex}{}
\titleformat{\subsection}[block]{\bfseries\filcenter}{\thesubsection.}{1ex}{}
\titleformat{\subsubsection}[runin]{\bfseries}{\thesubsubsection.}{1ex}{}[.]

\titlespacing*{\section}{0pt}{*3}{*1}
\titlespacing*{\subsection}{0pt}{*3}{*1}

\makeatletter
\newcommand{\neutralize}[1]{\expandafter\let\csname c@#1\endcsname\count@}
\makeatother

\newenvironment{theobis}[1]
{\renewcommand{\theequation}{\ref{#1}$'$}
	\neutralize{equation}\phantomsection
	\begin{theo}}
	{\end{theo}}

\renewcommand{\mkbibnamefirst}[1]{\textsc{#1}}
\renewcommand{\mkbibnamelast}[1]{\textsc{#1}}
\renewcommand{\mkbibnameprefix}[1]{\textsc{#1}}
\renewcommand{\mkbibnameaffix}[1]{\textsc{#1}}
\renewcommand{\finentrypunct}{}

\renewbibmacro{in:}{}

\renewbibmacro*{volume+number+eid}{%
	\printfield{volume}%
	\setunit*{\addnbspace}
	\printfield{number}%
	\setunit{\addcomma\space}%
	\printfield{eid}}

\DeclareFieldFormat[article]{volume}{\textbf{#1}\space}
\DeclareFieldFormat[article]{number}{\mkbibparens{#1}}

\DeclareFieldFormat{journaltitle}{#1,}
\DeclareFieldFormat[thesis]{title}{\mkbibemph{#1}\addperiod}
\DeclareFieldFormat[article, unpublished, thesis]{title}{\mkbibemph{#1},}
\DeclareFieldFormat[book]{title}{\mkbibemph{#1}\addperiod}
\DeclareFieldFormat[unpublished]{howpublished}{#1, }

\DeclareFieldFormat{pages}{#1}

\DeclareFieldFormat[article]{series}{Ser.~#1\addcomma}

\renewcommand*{\bibfont}{\small}

\begin{document}
	\pagestyle{plain}
	
	\maketitle
	
	\begin{abstract}
		Let $\G$ be a countably infinite group. A common theme in ergodic theory is to start with a probability measure\-/preserving \ep{\pmp} action $\G \acts (X, \mu)$ and a map $f \in L^1(X, \mu)$, and to compare the {global average} $\int f \,\Diff\mu$ of $f$ to the {pointwise averages} $|D|^{-1} \sum_{\delta \in D} f(\delta \cdot x)$, where $x \in X$ and $D$ is a nonempty finite subset of $\G$. The basic hope is that, when $D$ runs over a suitably chosen infinite sequence, these pointwise averages should converge to the global value for $\mu$-almost all $x$.
		
		In this paper we prove several results that refine the above basic paradigm by uniformly controlling the averages over {specific} sets $D$ rather than considering their limit as $|D| \to \infty$. Our results include ergodic theorems for the Bernoulli shift action $\G \acts ([0;1]^\G, \lambda^\G)$ and strengthenings of the theorem of Ab\'ert and Weiss that the shift is weakly contained in every free \pmp action of $\G$. In particular, we establish a purely Borel version of the Ab\'ert--Weiss theorem for finitely generated groups of subexponential growth. The central role in our arguments is played by the recently introduced measurable versions of the Lov\'asz Local Lemma, due to the current author and to Cs\'oka, Grabowski, M\'ath\'e, Pikhurko, and Tyros.
	\end{abstract}
	
	\section{Introduction}
	
	The \emph{Lov\'asz Local Lemma} \ep{the \emph{LLL} for short} is a powerful tool in probabilistic combinatorics, introduced by Erd\H os and Lov\'asz \cite{EL}. The LLL is mostly used to obtain existence results, and it is particularly well-suited for showing that a given structure $X$ admits a coloring satisfying some ``local'' constraints. Roughly speaking, in order for the LLL to apply in this context, two requirements must be met: First, a \emph{random} coloring should be ``likely'' to fulfill each individual constraint; second, the constraints must not interact with each other ``too much.'' For the precise statement, see \S\ref{subsec:SLLL}.
	
	It has been a matter of interest to determine if the LLL can be used to derive conclusions that are, in some sense, ``constructive'' \ep{as opposed to pure existence results}. A decisive breakthrough was made by Moser and Tardos \cite{MT}, who developed an \emph{algorithmic} approach to the LLL. \ep{The work of Moser and Tardos was preceded by a line or earlier results, starting with Beck's paper \cite{Beck}; for more details, see the references in \cite{MT}.} The Moser--Tardos method proved quite versatile and was adapted to establish ``constructive'' analogs of the LLL in a variety of different contexts. For example, Rumyantsev and Shen \cite{RSh} proved a \emph{computable} version of the LLL. Here we will be focused on the \emph{measurable} versions of the LLL that were studied in \cite{MLLL} by the current author and in \cite{CGMPT} by Cs\'{o}ka, Grabowski, M\'ath\'e, Pikhurko, and Tyros \ep{see also \cite{Kun} for related work by Kun}.
	
	Measurable analogs of the LLL are designed to apply in the following framework. Let $(X, \mu)$ be a standard probability space and let $C$ be a set of colors \ep{we will only consider the case when $C$ is finite}. Suppose we are looking for a coloring $f \colon X \to C$ that fulfills a family $\B$ of constraints. Under suitable assumptions, the ordinary LLL implies that such a coloring $f$ exists; however, this $f$ need not behave well with respect to the measurable structure on $(X, \mu)$. In contrast to that, measurable versions of the LLL can provide a \emph{$\mu$-measurable} \ep{or sometimes even \emph{Borel}} function $f \colon X \to C$ that satisfies the constraints $\B$, or at least does so on a ``large'' subset of $X$. Such results appear to be particularly relevant in ergodic theory, since many concepts pertaining to measure\-/preserving group actions are phrased in terms of measurable partitions of the underlying probability space---which can naturally be thought of as measurable colorings. Some ergodic\-/theoretic applications of the LLL can be found in \cite{MLLL, LSS}. Here we present further consequences of the LLL in measurable dynamics, specifically in the study of ergodic averages and of weak containment of measure\-/preserving group actions.
	
	Our arguments employ a general approach that is standard in combinatorics, in particular in graph coloring theory \ep{see, e.g., the book \cite{MR} for many examples}. The first step is to use \emph{concentration of measure} to obtain strong upper bounds on probabilities of certain ``bad'' random events; the LLL is then invoked to eliminate all the ``bad'' events. Nontrivial results can also be derived by combining the concentration of measure bounds with more classical tools, such as the Borel--Cantelli lemma \ep{Theorem~\ref{theo:weak_erg} below is as an example}. Roughly speaking, using the LLL instead of the Borel--Cantelli lemma results in replacing pointwise convergence with approximation in the $\infty$-norm.
	
	\subsubsection*{{Acknowledgement}}
	
	I am very grateful to Anush Tserunyan for many insightful discussions and to the anonymous referee for helpful suggestions.
	
	\section{Statements of results}
	
	Throughout, $\G$ denotes a countably infinite group with identity element $\mathbf{1}$. We study \emphd{probability measure\-/preserving \ep{\pmp} actions} of $\G$, i.e., actions of the form $\alpha \colon \G \acts (X, \mu)$, where $(X, \mu)$ is a standard probability space and the measure $\mu$ is $\alpha$-invariant. We also consider, more generally, \emphd{Borel actions} $\alpha \colon \G \acts X$, i.e., actions of $\G$ on a standard Borel space $X$ by Borel automorphisms.
	
	Given a set $A$, the \emphd{shift action}
	$\sigma_A \colon \G \acts A^\G$ on the set of all maps $x \colon \G \to A$
	is defined by
	\[
	(\gamma \cdot x)(\delta) \defeq x(\delta \gamma) \qquad \text{for all } x \in A^\G \text{ and } \gamma,\ \delta \in \G.	
	\]
	We are particularly interested in the case when $A$ is the unit interval $[0;1]$ equipped with the Lebesgue probability measure $\lambda$. \ep{Owing to the measure isomorphism theorem \cite[Theorem~17.41]{K_DST}, any other atomless standard probability space could be used instead.}
	To unclutter the notation, set \[(\Omega, \LG) \defeq ([0;1]^\G, \lambda^\G)\] and $\sigma \defeq \sigma_{[0;1]}$.
	Note that the action $\sigma \colon \G \acts (\Omega, \LG)$ is measure\-/preserving.
	
	\subsection{Ergodic theorems for the shift action}\label{subsec:erg}
	
	Let $\alpha \colon \G \acts (X, \mu)$ be a \pmp action. Given $f \in L^1(X, \mu)$, we can compute its \emphd{global average}:
	\begin{equation*}
		\E_\mu f \defeq \int_X f \,\Diff \mu,
	\end{equation*}
	and compare it to the \emphd{pointwise averages} of the form
	\begin{equation*}
		\E_D f(x) \defeq \frac{1}{|D|} \sum_{\delta \in D} f(\delta \cdot x),
	\end{equation*}
	where $x \in X$ and $D$ is a nonempty finite subset of $\G$. Note that $\E_D \colon L^1(X, \mu) \to L^1(X, \mu)$ is a linear operator of norm $1$: The lower bound on $\|\E_D\|_{\mathrm{op}}$ is witnessed by the constant $1$ function, while the upper bound follows from the fact that, since $\mu$ is $\alpha$-invariant, we have $\E_\mu \E_D f = \E_\mu f$, and hence
	\[
		\|\E_D f\|_1 \,=\, \E_\mu |\E_D f| \,\leq\, \E_\mu \E_D |f| \,=\, \E_\mu |f| \,=\, \|f\|_1.
	\]
	
	Assuming the action $\alpha$ is ergodic, one hopes to show that the pointwise averages $\E_D f$ converge, in a suitable sense, to $\E_\mu f$, as $D$ ranges over a given infinite family of finite subsets of $\G$. Results of this kind are usually referred to as \emph{ergodic theorems} \ep{often with adjectives indicating the mode of convergence, such as ``pointwise''}. Two prototypical examples are \emph{von Neumann's} \cite{vN} and \emph{Birkhoff's} \cite{Birkhoff} \emph{ergodic theorems}. Both of these classical results apply when $\G = \Z$ and $D$ ranges over the sets of the form $\set{0, 1, \ldots, n-1}$ with $n \in \N^+$. Von Neumann's theorem yields convergence in the $2$-norm \ep{assuming $f \in L^2(X, \mu)$ to begin with}, while Birkhoff's result ensures pointwise convergence almost everywhere. An extension of Birkhoff's pointwise ergodic theorem to all amenable $\G$ was obtained by Lindenstrauss \cite{Lin}; there $D$ ranges over a \emph{tempered F\o lner sequence} \ep{the special case of Lindenstrauss's result for $f \in L^2(X, \mu)$ follows from the earlier work of Shulman, see \cite[\S5.6]{Temp}}. Generalizing ergodic theorems beyond the realm of amenable groups is a major challenge; for further background,
	see, e.g., \cite{AAB, BufKlim, BowNev} and the references therein.
	
	Here we work with an \emph{arbitrary} group $\G$; moreover, the only condition on the sequence $(D_n)_{n \in \N}$ of averaging sets is that $|D_n|$ grows sufficiently quickly with $n$. On the other hand, instead of studying arbitrary ergodic actions, we focus our attention on the shift action $\sigma \colon \G \acts (\Omega, \LG)$ in the hope of exploiting its mixing properties. Our first result is a pointwise ergodic theorem for \emph{continuous} functions $f \colon \Omega \to \C$:
	
	\begin{theo}[\textls{Pointwise ergodic theorem for continous maps on the shift}]\label{theo:weak_erg}
		Let $(D_n)_{n \in \N}$ be a sequence of finite subsets of $\G$ such that $|D_n|/\log n \to\infty$. Then, for all continuous $f \colon \Omega \to \C$,
		\[
		\lim_{n \to \infty} \E_{D_n} f(x) \,=\, \E_{\LG} f, \qquad \text{for $\LG$-a.e.} \ x \in \Omega.
		\]
	\end{theo}
	
	Since the set of all continuous functions is dense in $L^1(\Omega, \LG)$ and $\|\E_D\|_{\mathrm{op}} = 1$ for all nonempty finite $D \subset \G$, Theorem~\ref{theo:weak_erg} has the following immediate corollary:
	
	\begin{corl}[\textls{Mean ergodic theorem for the shift}]\label{corl:mean_erg}
		Let $(D_n)_{n \in \N}$ be a sequence of finite subsets of $\G$ such that $|D_n|/\log n \to\infty$. Then, for all $f \in L^1 (\Omega, \LG)$, we have
		\[
		\lim_{n \to \infty} \E_{D_n} f \,=\, \E_{\LG} f \qquad \text{in}\  L^1(\Omega, \LG).
		\]
	\end{corl}
	
	It is natural to ask whether Theorem~\ref{theo:weak_erg} can be extended to all $f \in L^1(\Omega, \LG)$. The answer turns out to be negative even if the lower bound on the growth rate of the averaging sets is raised, as the constructions of Akcoglu and del~Junco \cite{AkcogluJunco} and del~Junco and Rosenblatt \cite{JunRos} \ep{with minor modifications} demonstrate:
	
	\begin{theo}[{$\text{ess.}$ Akcoglu--del Junco \cite{AkcogluJunco} and del Junco--Rosenblatt \cite{JunRos}}]\label{theo:bad}
		Suppose that $\G = \Z$ and let $h \colon \N \to \N$ be an arbitrary function. There exists a sequence $(D_n)_{n \in \N}$ of finite subsets of $\Z$ with the following properties:
		\begin{itemize}[label=--,wide]
			\item each $D_n$ is an interval, i.e., a set of the form $\set{s, s+1, \ldots, s + \ell - 1}$ for $s \in \Z$ and $\ell \in \N^+$;
			
			\item $|D_n| \geq h(n)$ for all $n \in \N$;
			
			\item for every free \pmp action $\Z \acts (X, \mu)$, there is a Borel set $A \subseteq X$ such that
			\[
				\liminf_{n \to \infty} \,\E_{D_n} \mathbbm{1}_A(x) \,=\, 0 \quad \text{and} \quad \limsup_{n \to \infty} \, \E_{D_n} \mathbbm{1}_A(x) \,=\, 1, \qquad \text{for $\mu$-a.e.}\ x \in X,
			\]
			where $\mathbbm{1}_A \colon X \to \set{0,1}$ is the indicator function of $A$. Moreover, the family of such sets $A$ is comeager in the measure algebra $\operatorname{MAlg}(X, \mu)$.
		\end{itemize}
	\end{theo}
	
	For completeness, we sketch a proof of Theorem~\ref{theo:bad} using Rokhlin's lemma in the \hyperref[sec:app]{appendix}.
	
	As mentioned in the introduction, Theorem~\ref{theo:weak_erg} follows by combining a concentration of measure inequality with the Borel--Cantelli lemma. We now turn to further results that can be obtained if the Borel--Cantelli lemma is replaced by the LLL.
	
	For a \pmp action $\alpha \colon \G \acts (X, \mu)$, $f \in L^1(X, \mu)$, and a nonempty finite set $D \subset \G$, define the \emphd{discrepancy norm} of $f$ with respect to $D$ by the formula
	\[
		\|f\|^{\mathrm{disc}}_D \defeq \| \E_D f \,-\, \E_\mu f\|_\infty.
	\]
	\ep{Here $\| \cdot \|_\infty$ is the $\infty$-norm in the sense of $L^\infty(X, \mu)$.} Even if $f \colon \Omega \to \C$ is continuous, its discrepancy norm may be separated from $0$.
	For instance, consider the continuous map \[f \colon \Omega \to [-1;1] \colon x \mapsto -1 + 2\cdot x(\mathbf{1}).\]
	Then $\E_\LG f = 0$ yet $\| \E_D f\|_\infty = 1$, and hence $\|f\|^{\mathrm{disc}}_D = 1$, for all nonempty finite $D \subset \G$. However, we show that any $f \in L^1(\Omega, \LG)$ can be written as a {sum} of two functions $g$, $h \in L^1(\Omega, \LG)$, where $g$ is small in the discrepancy norm, while $h$ is small in the $1$-norm:
	
	\begin{theo}[\textls{$L^\infty$-ergodic theorem for the shift}]\label{theo:strong_erg}
		For all $f \in L^1(\Omega, \LG)$ and $\epsilon > 0$, there exists $C > 0$ with the following property:
		
		\smallskip
		
		Let $(D_n)_{n \in \N}$ be a sequence of finite subsets of $\G$ with $|D_n|\geq C \log(n+2)$ for all $n \in \N$. Then there exist $g$, $h \in L^1(\Omega, \LG)$ such that $f = g + h$, $\|h\|_1 \leq \epsilon$, and $\|g\|^{\mathrm{disc}}_{D_n} \leq \epsilon$ for all $n \in \N$.
	\end{theo}
	
	Note that Theorem~\ref{theo:strong_erg} also yields Corollary~\ref{corl:mean_erg}.
	
	Our ultimate goal in this subsection is to sharpen Theorem~\ref{theo:strong_erg} by considering other statistical properties of the function $f$, beside its average $\E_\LG f$. This is made precise by the following formalism. Let $K$ be a compact metric space. We use $\Prob(K)$ to denote the space of all probability Borel measures on $K$ equipped with the usual weak-$\ast$ topology \ep{see, e.g., \cite[\S{}17.E]{K_DST}}. Given a \pmp action $\alpha \colon \G \acts (X, \mu)$ and a Borel function $f \colon X \to K$, define
	\[
		\M_\mu f \defeq f_\ast(\mu),
	\]
	where $f_\ast \colon \Prob(X) \to \Prob(K)$ is the pushforward map. \ep{This notation is intended to be reminiscent of $\E_\mu f$, while the letter ``$\M$'' emphasizes that $\M_\mu f$ is a \emph{measure}.} For $x \in X$ and a nonempty finite set $D \subset \G$, let $\upsilon_{x, D}$ be the probability measure on $X$ with \ep{finite} support $D \cdot x$ given by
	\begin{equation*}
	\uni_{x, D}(\set{y}) \defeq \frac{1}{|D|} \cdot |\set{\delta \in D \,:\, \delta \cdot x = y}| \qquad \text{for all } y \in D \cdot x.
	\end{equation*}
	If the $\alpha$-stabilizer of $x$ is trivial, then $\uni_{x, D}$ is simply the uniform probability measure on $D \cdot x$. Let
	\[
		\M_D f(x) \defeq f_\ast(\uni_{x, D}).
	\]
	The measures $\M_\mu f$ and $\M_D f(x)$ are points in $\Prob(K)$ that encode the global and the pointwise statistics of $f$, respectively. In particular, when $K$ is a subset of $\C$, $\M_\mu f$ and $\M_D f(x)$ contain the information about $\E_\mu f$ and $\E_D f(x)$; explicitly,
	\begin{equation}\label{eq:star}
		\E_\mu f = \int_K z \,\Diff (\M_\mu f)(z) \qquad \text{and} \qquad \E_D f(x) = \int_K z \,\Diff (\M_D f(x))(z).
	\end{equation}
	
	We wish to also take into account more detailed information about the interaction of $f$ with the action $\alpha$. Toward that end, let $\code{f} \colon X \to K^\G$ denote the equivariant map given by
	\[
		\code{f}(x)(\gamma) \defeq f(\gamma \cdot x) \qquad \text{for all } x \in X \text{ and } \gamma \in \G.
	\]
	The map $\code{f}$ is called the \emphd{symbolic representation}, or the \emphd{coding map}, of the dynamical system $(X, \G, \alpha, f)$. Notice that the projection function $\pr \colon K^\G \to K \colon \kappa \mapsto \kappa(\mathbf{1})$ satisfies $f = \pr \circ \code{f}$ and gives rise to a continuous map $\pr_\ast \colon \Prob(K^\G) \to \Prob(K)$ such that $f_\ast = \pr_\ast \circ (\code{f})_\ast$.
	This observation shows that, by considering $\code{f}$, we achieve greater generality than just by working with $f$ itself.
	
	Given a standard probability space $(X, \mu)$ and a compact metric space $(K, \dist)$, let $\mathfrak{B}(X, K)$ denote the set of all Borel functions from $X$ to $K$. We equip $\mathfrak{B}(X, K)$ with a psedometric ${\dist}_\mu$ given by
	\[
		{\dist}_\mu(f, g) \defeq \int_X \dist(f(x), g(x))\, \Diff\mu(x).
	\]
	If $K$ is a subset of $\C$ equipped with the metric $\dist(z_1, z_2) = |z_1 - z_2|$, then ${\dist}_\mu(f, g) = \|f - g\|_1$.
	
	\begin{theo}[\textls{Pushforward-ergodic theorem for the shift}]\label{theo:ult_erg}
		Let $(K, \dist)$ be a compact metric space and let $f \colon \Omega \to K$ be a Borel function. For any $\epsilon > 0$ and an open neighborhood $U$ of the measure $\M_\LG\code{f}$, there exists $C > 0$ with the following property:
		
		\smallskip
		
		Let $(D_n)_{n \in \N}$ be a sequence of finite subsets of $\G$ with $|D_n|\geq C \log(n+2)$ for all $n \in \N$.
		Then there is a Borel map $g \colon \Omega \to K$ such that ${\dist}_\LG(f, g) \leq \epsilon$ and
		\[
		\M_{D_n} \code{g}(x) \in U, \qquad \text{for all}\ n \in \N\ \text{and for $\LG$-a.e.}\ x \in X.
		\]
	\end{theo}
	
	In the light of \eqref{eq:star}, it is clear that Theorem~\ref{theo:strong_erg} is a special case of Theorem~\ref{theo:ult_erg}. 
	
	We end this subsection with a simple application of Theorem~\ref{theo:ult_erg}. Recall that a group $\G$ is called \emphd{residually finite} if the intersection of all its subgroups of finite index is trivial. The following is an easy observation:
	
	\begin{prop}\label{prop:res_fin}
		A countable group $\G$ is residually finite if and only if every open neighborhood $U \subseteq \Prob(\Omega)$ of $\LG$ contains a finitely supported measure $\nu$ that is shift-invariant.
	\end{prop}
	\begin{coolproof}
		Let $(\nu_n)_{n \in \N}$ be a sequence of finitely supported shift-invariant measures on $\Omega$ that converges to $\LG$. This gives us a sequence of actions of $\G$ on the finite sets $X_n \defeq \supp(\nu_n)$ and, since $\nu_n \to \LG$, each nonidentity group element $\gamma \in \G$ acts on $X_n$ nontrivially for all large enough $n$. This shows that $\G$ is residually finite.
		
		Conversely, suppose that $\G$ is residually finite and let $(\Delta_n)_{n \in \N}$ be a decreasing sequence of finite index subgroups of $\G$ with trivial intersection. For $k \in \N^+$, let
		\[
			Q_k \defeq \set{0,\, 1/k,\, \ldots,\, (k-1)/k} \subset [0;1].
		\]
		Let $P(k, n)$ denote the set of all maps $x \colon \G \to Q_k$ that are constant on the right cosets of $\Delta_n$. Then the set $P(k, n)$ is finite and shift-invariant, and, letting $\nu_{k, n}$ be the uniform probability measure on $P(k,n)$, we see that $\nu_{k, n} \to \LG$ as $k$, $n \to \infty$.
 	\end{coolproof}
	
	Motivated by Proposition~\ref{prop:res_fin}, we say that a group $\G$ is \emphd{approximately residually finite} if for every open neighborhood $U$ of $\LG$, there is a finitely supported measure $\nu$ such that $\gamma \cdot \nu \in U$ for all $\gamma \in \G$. Proposition~\ref{prop:res_fin} implies that every residually finite group is approximately residually finite, so our terminology is consistent. An intuitive way of thinking about approximate residual finiteness is as follows: To show that a group $\G$ is approximately residually finite, we have to find finite subsets $X \subset \Omega$ that are ``almost uniformly distributed'' over the space $(\Omega, \LG)$ and also remain such when shifted by any $\gamma \in \G$. We remark that a \emph{random} finite set $X$ fails to have this property: For any $n \in \N^+$, the product action $\sigma^n \colon \G \acts (\Omega^n, \LG^n)$ is ergodic, and hence if $x_1$, \ldots, $x_n \in \Omega$ are chosen randomly and independently from each other, then, with probability $1$, for every open $V \subseteq \Omega$ there is some $\gamma \in \G$ such that $\gamma \cdot x_1$, \ldots, $\gamma \cdot x_n \in V$.  Nevertheless, we have the following:
	
	\begin{corl}[to Theorem~\ref{theo:ult_erg}]\label{corl:approx_res_fin}
		Every countable group is approximately residually finite.
	\end{corl}
	\begin{coolproof}
		Let $U$ be an open neighborhood of $\LG$. It suffices to exhibit a finitely supported measure $\nu \in U$ such that $\nu \cdot \gamma \in U$ for all $\gamma \in \G$, where the \emphd{right shift action} $\Omega \curvearrowleft \G$ is given by
		\[
		(x \cdot \gamma)(\delta) \defeq x(\gamma\delta) \qquad \text{for all } x \in \Omega \text{ and } \gamma,\ \delta \in \G.	
		\]
		Applying Theorem~\ref{theo:ult_erg} with $K = [0;1]$ and $f = (x \mapsto x(\mathbf{1}))$, we obtain a nonempty finite set $D \subset \G$ and a Borel map $g \colon \Omega \to [0;1]$ such that $\M_D \code{g}(x) \in U$ for $\LG$-a.e.\ $x \in \Omega$. Since $\G$ is countable and the measure $\LG$ is right-shift-invariant, there is $x \in \Omega$ such that $\M_D \code{g}(x \cdot \gamma) \in U$ for all $\gamma \in \G$. Set $\nu \defeq \M_D \code{g}(x)$. Then $\nu$ is finitely supported; furthermore, it is straightforward to verify, using the \ep{left-}equivariance of $\code{g}$ and the fact that the left and the right shift actions of $\G$ on $\Omega$ commute with each other, that $\nu \cdot \gamma = \M_D\code{g}(x \cdot \gamma)$ for all $\gamma \in \G$. Hence, $\nu$ is as desired.
	\end{coolproof}
	
	Since the above argument only involves the properties of $g$ on a countable subset of $\Omega$, Corollary~\ref{corl:approx_res_fin} can also be derived directly from the classical LLL, without using its measurable analogs.
	
	\subsection{Pointwise versions of the Ab\'ert--Weiss theorem}\label{subsec:AW}
	
	So far we have considered the action $\sigma \colon \G \acts (\Omega, \LG)$ on its own. Now we would like to discuss the relationship between $\sigma$ and other actions of $\G$.
	
	The concepts of \emph{weak containment} and \emph{weak equivalence} of \pmp actions were introduced by Kechris in~\cite[\S10(C)]{K_book}. They are inspired by the analogous notions for unitary representations and are closely related to the so-called \emph{local\-/global convergence} in the theory of graph limits~\cite{LocalGlobal}. The relation of weak equivalence is much coarser than the isomorphism relation, which makes it relatively well-behaved. On the other hand, several interesting parameters associated with \pmp actions---such as their cost, type, etc.---turn out to be invariants of weak equivalence. Due to these favorable properties, the relations of weak containment and weak equivalence have attracted a considerable amount of attention in recent years. For a survey of the topic, see \cite{BK}.
	
	Roughly speaking, a \pmp action $\alpha \colon \G \acts (X, \mu)$ is weakly contained in another \pmp action $\beta \colon \G \acts (Y, \nu)$ if for every compact metric space $K$ and for any Borel map $f \colon X \to K$, the interaction of $f$ with $\alpha$ can be arbitrarily well ``simulated'' by a Borel map $g \colon Y \to K$ interacting with $\beta$. Here is a precise definition:
	
	\begin{defn}[{\textls{Weak containment}; {\cite[\S{}2.2(2)]{BK}}}]
		Let $\alpha \colon \G \acts (X, \mu)$ and $\beta \colon \G \acts (Y, \nu)$ be \pmp actions of $\G$. We say that $\alpha$ is \emphd{weakly contained} in $\beta$, in symbols $\alpha \preceq \beta$, if for any compact metric space $K$, a Borel function $f \colon X \to K$, and an open neighborhood $U$ of the measure $\M_\mu \code{f}$, there exists a Borel map $g \colon Y \to K$ such that $\M_\nu \code{g} \in U$. If both $\alpha \preceq \beta$ and $\beta \preceq \alpha$, then $\alpha$ and $\beta$ are said to be \emphd{weakly equivalent}, in symbols $\alpha \simeq \beta$.
	\end{defn}

	Weak containment can be defined in a number of equivalent ways, several of which can be found in~\cite[\S\S2.1, 2.2]{BK}. The characterization given above is due to Ab\'ert and Weiss~\cite[Lemma~8]{AW} \ep{see also \cite[Proposition~3.5]{T-D}}. We sometimes write $(\alpha, \mu) \preceq (\beta, \nu)$ instead of $\alpha \preceq \beta$ in order to emphasize the dependence of weak containment on the invariant measures $\mu$ and $\nu$.
	
	Burton~\cite[Corollary~4.2]{B} (see~\cite[Theorem~3.3]{BK}) proved that there exist continuum many distinct weak equivalence classes of \ep{not necessarily ergodic} \pmp actions of $\G$.
	Glasner, Thouvenot, and Weiss~\cite{GTW} and independently Greg Hjorth (unpublished) proved that the pre-order of weak containment has a \emph{maximum} element (see also~\cite[Theorem~3.1]{BK}). A complementary result of Ab\'ert and Weiss~\cite[Theorem~1]{AW} (see also~\cite[Theorem~3.5]{BK}) asserts that the shift action $\sigma \colon \G \acts (\Omega, \LG)$ is \emph{minimum} among all \pmp actions $\alpha \colon \G \acts (X, \mu)$ that are \emphd{\ep{almost everywhere} free}, i.e.,  such that the $\alpha$-stabilizer of $\mu$-a.e.\ $x \in X$ is trivial:

	\begin{theo}[{Ab\'ert--Weiss~{\normalfont\cite[Theorem~1]{AW}}}]\label{theo:AW}
		Let $\alpha \colon \G \acts (X, \mu)$ be an almost everywhere free \pmp action of $\G$. Then $(\sigma, \LG) \preceq (\alpha, \mu)$; or, explicitly, the following statement holds:
		
		\smallskip
		
		Let $K$ be a compact metric space and let $f \colon \Omega \to K$ be a Borel function. Then, for any open neighborhood $U$ of the measure $\M_\LG \code{f}$, there is a Borel map $g \colon X \to K$ such that $\M_\mu \code{g} \in U$.
	\end{theo}
	
	We strengthen Theorem~\ref{theo:AW} by replacing the measure $\M_\mu \code{g}$ by its pointwise analogs of the form $\M_D \code{g} (x)$. Moreover, our result applies to actions that are not necessarily free but only ``close enough'' to being free. Specifically, for a set $S \subseteq \G$, we say that an action $\alpha \colon \G \acts X$ is \emphd{$S$-free} if for all $\gamma$, $\delta \in S$ and $x \in X$,
	$\gamma \cdot x = \delta \cdot x$ implies $\gamma = \delta$. \ep{Thus, ``free'' is the same as ``$\G$-free.''} Given a sequence of sets $S_1$, \ldots, $S_n \subseteq \G$, we say that $\alpha$ is \emphd{$(S_1, \ldots, S_n)$-free} if $\alpha$ is $S_i$-free for each $1 \leq i \leq n$. 
	
	\begin{theo}[\textls{Pointwise Ab\'ert--Weiss}]\label{theo:main_delta}
		Let $K$ be a compact metric space and let $f \colon \Omega \to K$ be a Borel function. For any open neighborhood $U$ of the measure $\M_\LG \code{f}$, there exist $C > 0$ and a finite set $S \subset \G$ with the following property:
		
		\smallskip
		
		Let $D$ be a finite subset of $\G$ with $|D| \geq C$ and let $\alpha \colon \G \acts X$ be an $(S, D)$-free Borel action of $\G$. Then, for any $\mu \in \Prob(X)$ and $\delta > 0$, there is a Borel map $g \colon X \to K$ such that
		\[
		\mu(\set{x \in X \,:\, \M_D \code{g} (x) \in U}) \,\geq\, 1 - \delta.
		\]
	\end{theo}

	\begin{remks}\label{remk:delta}
		Let us make a few comments about the statement of Theorem~\ref{theo:main_delta}.
		
		\begin{enumerate}[label={\ep{\roman*}},wide]
			\item To see that Theorem~\ref{theo:main_delta} is a strengthening of the \hyperref[theo:AW]{Ab\'ert--Weiss theorem}, let $\alpha \colon \G \acts (X, \mu)$ be a free \pmp action. Given a compact metric space $K$ and a Borel function $f \colon \Omega \to K$, we can apply Theorem~\ref{theo:main_delta} to obtain a finite set $D \subset \G$ and a Borel map $g \colon X \to K$ such that the pushforward measure $\M_D \code{g} (x)$ is arbitrarily close to $\M_\LG \code{f}$, for all points $x \in X$ away from a set of arbitrarily small measure. The $\alpha$-invariance of $\mu$ yields
			\[
				\mu = \int_X \uni_{x,D} \,\Diff \mu(x), \quad \text{hence} \quad \M_\mu \code{g} = \int_X \M_D \code{g} (x) \,\Diff \mu(x), 
			\]
			and thus $\M_\mu \code{g}$ is also close to $\M_\LG \code{f}$, as desired.
			
			\smallskip
			
			\item The measure $\mu$ in Theorem~\ref{theo:main_delta} is not required to be $\alpha$-invariant \ep{or even $\alpha$-quasi-invariant} and is only used to bound the set of all $x \in X$ with $\M_D \code{g} (x) \not \in U$.
			
			\smallskip
			
			\item We emphasize that the averaging set $D$ in Theorem~\ref{theo:main_delta} is independent of the choice of $\delta > 0$; that is what makes this result particularly interesting. It is possible that the conclusion of Theorem~\ref{theo:main_delta} also holds with $\delta = 0$, but we do not know how to prove \ep{or disprove} that in general; see Problem~\ref{prob:delta_zero} in Section \ref{sec:remks}. \ep{However, we can make $\delta$ be zero under some additional assumptions---see \ref{item:shift} and Theorem~\ref{theo:main_Borel} below.}
			
			\smallskip
			
			\item\label{item:shift} In contrast to the \hyperref[theo:AW]{Ab\'ert--Weiss theorem}, the conclusion of Theorem~\ref{theo:main_delta} is nontrivial even if $(X, \mu) = (\Omega, \LG)$ and $\alpha = \sigma$. This case, however, is already covered by the ergodic Theorem~\ref{theo:ult_erg}, in fact even with $\delta = 0$.
			
			\smallskip
			
			\item\label{item:derivation} For actions $\alpha$ that are free and measure\-/preserving, Theorem~\ref{theo:main_delta} follows relatively straightforwardly by combining Theorem~\ref{theo:ult_erg} with the usual \hyperref[theo:AW]{Ab\'ert--Weiss theorem}. We sketch the argument here. Let $\alpha \colon \G \acts (X, \mu)$ be a free \pmp action. Let $K$ be a compact metric space and let $f \colon \Omega \to K$ be a Borel function. Fix an open neighborhood $U$ of the measure $\M_\LG \code{f}$. By Theorem~\ref{theo:ult_erg}, for any sufficiently large finite set $D \subset \G$, there is a Borel map $h \colon \Omega \to K$ with
			\begin{equation}\label{eq:h}
				\M_D \code{h} (x) \in U, \qquad \text{for $\LG$-a.e.}\ x \in \Omega.
			\end{equation}
			The equivariance of $\code{h}$ yields $\M_D \code{h}(x) = (\code{h})_\ast(\uni_{x, D}) = \uni_{\code{h}(x), D}$, and hence \eqref{eq:h} is equivalent to
			\[
				\uni_{\kappa, D} \in U, \qquad \text{for $\M_\LG \code{h}$-a.e.}\ \kappa \in K^\G. 
			\]
			Now we can use the \hyperref[theo:AW]{Ab\'ert--Weiss theorem} to obtain a Borel map $g \colon X \to K$ for which the pushforward measure $\M_\mu \code{g}$ is so close to $\M_\LG \code{h}$ that
			\[
				\mu(\set{x \in X \,:\, \M_D \code{g} (x) \in U}) \,=\, \M_\mu \code{g} (\set{\kappa \in K^\G \,:\, \uni_{\kappa, D} \in U}) \,\geq\, 1-\delta,
			\]
			for any given $\delta > 0$, as desired. For non-free actions $\alpha$, a different, more direct proof is necessary.
			
			\smallskip
			
			\item The results of \S\ref{subsec:erg} apply to an infinite sequence of averaging sets $(D_n)_{n \in \N}$, while in Theorem~\ref{theo:main_delta} we only consider a single set $D$. Our approach can be routinely adapted to extend Theorem~\ref{theo:main_delta} to the case of finitely many averaging sets; however, when the family of averaging sets is infinite, our methods are not applicable---see Remark \ref{remk:amenable}.
		\end{enumerate}
	\end{remks}	
	
	Notice that the pointwise operator $\M_D$ is well-defined for an arbitrary Borel action $\alpha \colon \G \acts X$ and does not require fixing a probability measure $\mu$ on $X$. Therefore, it makes sense to ask for a \emph{purely Borel} version of the \hyperref[theo:AW]{Ab\'ert--Weiss theorem}, with the last line of Theorem~\ref{theo:main_delta} replaced by
	\[
		\text{$\M_D \code{g} (x) \in U$, \qquad \textsl{for all} $x \in X$.}
	\]
	Here we establish such a version for finitely generated groups of subexponential growth and, more generally, for uniformly subexponential Borel actions. Let $\alpha \colon \G \acts X$ be a Borel action of $\G$. We say that $\alpha$ is \emphd{uniformly subexponential} if for every finite set $S \subset \G$ and for all $\epsilon > 0$, there is $n_0 \in \N$ such that for all $n \geq n_0$ and for all $x \in X$, $\left|S^n \cdot x\right| \leq (1 + \epsilon)^n$, where $S^n \defeq \set{\gamma_1 \cdots \gamma_n \,:\, \gamma_i \in S \text{ for all }1 \leq i \leq n}$. For example, if $\G$ is a finitely generated group of subexponential growth, then every action of $\G$ is uniformly subexponential.
	
	\begin{theo}[\textls{Borel Ab\'ert--Weiss for uniformly subexponential actions}]\label{theo:main_Borel}
		Let $K$ be a compact metric space and let $f \colon \Omega \to K$ be a Borel function. For any open neighborhood $U$ of the measure $\M_\LG \code{f}$, there exist $C > 0$ and a finite set $S \subset \G$ with the following property:
		
		\smallskip
		
		Let $D$ be a finite subset of $\G$ with $|D| \geq C$ and let $\alpha \colon \G \acts X$ be a uniformly subexponential $(S, D)$-free Borel action of $\G$. Then there is a Borel map $g \colon X \to K$ such that \[\M_D \code{g} (x) \in U, \qquad \text{for all}\ x \in X.\]
	\end{theo}
	
	Note that, even though groups of subexponential growth are amenable, the averaging set $D$ in the statement of Theorem~\ref{theo:main_Borel} is {not} assumed to be a F\o{}lner set. 
	
	\subsection{Outline of the remainder of the paper}
	
	This paper is organized as follows. Section~\ref{sec:prelim} contains a few definitions and some preliminary results concerning the continuity of various basic operations, such as $f \mapsto f_\ast$. We commence the proofs of Theorems~\ref{theo:weak_erg}, \ref{theo:ult_erg}, \ref{theo:main_delta}, and \ref{theo:main_Borel} in Section~\ref{sec:red}, where they are reduced to their special cases with a more ``combinatorial'' flavor. Then, in Section~\ref{sec:concentration}, we state and prove a certain concentration of measure inequality. At this point, we already have all the tools needed to derive Theorem~\ref{theo:weak_erg}, which is done in \S\ref{subsec:weak_erg}. In Section~\ref{sec:LLL} we review the LLL and its measurable analogs, and in Section~\ref{sec:proofs} we complete the proofs of Theorems \ref{theo:ult_erg}, \ref{theo:main_delta} and \ref{theo:main_Borel}. It turns out that in order to prove Theorem~\ref{theo:ult_erg}, it is not enough to simply apply a known measurable version of the LLL---we actually have to go through the \emph{proof} of one of them to obtain some additional information; this is done in \S\ref{subsec:ult_erg}. We conclude the paper with some open problems in Section \ref{sec:remks}. The \hyperref[sec:app]{appendix} contains a proof of Theorem~\ref{theo:bad}.

	\section{Preliminaries}\label{sec:prelim}
	
	\subsection{Further notation}
	
	\subsubsection*{Integers}
	We use $\N$ to denote the set of all nonnegative integers and identify each $k \in \N$ with the $k$-element set $\set{i \in \N \,:\, i < k}$. Let $\N^+ \defeq \N \setminus \set{0}$. All finite sets \ep{including each $k \in \N$} are assumed to carry discrete topologies.
	
	\subsubsection*{Sets and functions}
	Each function $f$ is identified with its graph, i.e., the set $\set{(x, y) \,:\, y = f(x)}$. This enables the use of set-theoretic notation, such as $\subseteq$, $|\cdot|$, etc., for functions. For a function $f$ and a set $S$ of its domain, $\rest{f}{S}$ denotes the restriction of $f$ to $S$.
	For sets $A$ and $B$,
	\[
	\begin{array}{rll}
	\text{--} & \fins{B} & \text{denotes the set of all finite subsets of $B$;}\\
	\text{--} & [B \to A] & \text{denotes the set of all partial functions $\phi \colon B \rightharpoonup A$;}\\
	\text{--} & \finf{B}{A} & \text{denotes the set of all partial functions $\phi \colon B \rightharpoonup A$ with $\dom(\phi) \in \fins{B}$.}
	\end{array}
	\]
	The identity function $X \to X$ on a set $X$ is denoted by $\id_X$.

	\subsubsection*{Symbolic dynamics}
	
	Let $A$ be a set and let $\alpha \colon \G \acts X$ be an action of $\G$. We extend the definition of the coding map to partial functions $f \colon X \rightharpoonup A$ by letting $\pi_f(x) \colon \G \rightharpoonup A$ be given by
	\[
		\pi_f(x)(\gamma) \defeq \begin{cases}
		f(\gamma \cdot x) &\text{if } \gamma \cdot x \in \dom(f);\\
		\text{undefined} &\text{otherwise},
		\end{cases} \qquad \text{for all}\ x \in X \text{ and } \gamma \in \G.
	\]
	We similarly extend the shift action $\sigma_A \colon \G \acts A^\G$ to an action $\G \acts [\G \to A]$ in the obvious way.
	
	\subsubsection*{The free part of an action}
	
	For an action $\alpha \colon \G \acts X$ of $\G$, let $\Free(X) \subseteq X$ denote the set of all $x \in X$ whose $\alpha$-stabilizer is trivial 
	and let $\Free(\alpha) \colon \G \acts \Free(X)$ denote the induced action of $\G$ on $\Free(X)$; we call $\Free(\alpha)$ the \emphd{free part} of $\alpha$.
	
	\subsubsection*{Miscellaneous}
	
	For a metric space $(K, \dist)$, $a$, $b \in K$, and $\epsilon > 0$, we write $a \approx_\epsilon b$ to mean $\dist(a,b) < \epsilon$.
	
	\subsection{Topological preliminaries}

	\subsubsection*{Continuity of the coding map}
	
	Fix an arbitrary enumeration $\set{\gamma_n}_{n \in \N}$ of the elements of $\G$. If $(K, \dist)$ is a compact metric space, then the product topology on $K^\G$ is induced by the metric $\hat{\dist}$:
	\[
		\hat{\dist}(\kappa, \eta) \defeq \sum_{n=0}^\infty \frac{\dist(\kappa(\gamma_n), \eta(\gamma_n))}{2^{n+1}}.
	\]
	Recall that if $X$ is a standard Borel space and $\mu \in \Prob(X)$, then the space $\mathfrak{B}(X, K)$ is endowed with the pseudometric $\dist_\mu$. Additionally, we shall consider the uniform metric $\dist_{\mathsf{uni}}$ given by
	\[
		\dist_{\mathsf{uni}} (f, g) \defeq \sup_{x \in X} \dist (f(x), g(x)).
	\]
	
	\begin{lemma}\label{lemma:cont_of_code}
		Let $(K, \dist)$ be a compact metric space and let $\alpha \colon \G \acts X$ be a Borel action of $\G$.
		\begin{enumerate}[label=\ep{\normalfont\alph*},wide]
			\item\label{item:mu} If $\mu \in \Prob(X)$ is $\alpha$-invariant, then the function
			$(\mathfrak{B}(X, K), \dist_\mu) \to (\mathfrak{B}(X, K^\G), \hat{\dist}_\mu) \colon f \mapsto \pi_f$ is distance\-/preserving, hence continuous.
			
			\smallskip
			
			\item\label{item:uni} The function
			$(\mathfrak{B}(X, K), \dist_{\mathsf{uni}}) \to (\mathfrak{B}(X, K^\G), \hat{\dist}_{\mathsf{uni}}) \colon f \mapsto \pi_f
			$ is $1$-Lipschitz, hence continuous.
		\end{enumerate}
	\end{lemma}
	\begin{coolproof}
		\ref{item:mu} For all $f$, $g \in \mathfrak{B}(X, K)$, we have
		\[
		\hat{\dist}_\mu(\code{f}, \code{g}) \,=\, \int_X \hat{\dist}(\code{f}(x), \code{g}(x)) \,\Diff\mu(x) \,=\, \int_X \sum_{n=0}^\infty \frac{\dist(f(\gamma_n \cdot x), g(\gamma_n \cdot x))}{2^{n+1}} \,\Diff\mu(x).
		\]
		Switching the order of integration and summation, we rewrite the last expression as
		\[
		\sum_{n=0}^\infty \frac{1}{2^{n+1}}\int_X \dist(f(\gamma_n \cdot x), g(\gamma_n \cdot x)) \,\Diff\mu(x).
		\]
		Since $\mu$ is $\alpha$-invariant, this is equal to
		\[
		\sum_{n=0}^\infty \frac{1}{2^{n+1}} \int_X \dist(f(x), g(x)) \,\Diff \mu(x) \,=\, \sum_{n=0}^\infty\frac{\dist_\mu(f, g)}{2^{n+1}} \,=\, \dist_\mu(f, g).
		\]
		
		\ref{item:uni} For all $f$, $g \in \mathfrak{B}(X, K)$ and $x \in X$, we have
		\[
			\hat{\dist}(\code{f}(x), \code{g}(x)) \,=\, \sum_{n = 0}^\infty \frac{\dist(f(\gamma_n \cdot x), g(\gamma_n \cdot x))}{2^{n+1}} \,\leq\, \sum_{n = 0}^\infty \frac{\dist_{\mathsf{uni}}(f, g)}{2^{n+1}} \,=\, \dist_{\mathsf{uni}}(f, g),
		\]
		and the desired conclusion follows.
	\end{coolproof}
	
	\subsubsection*{Continuity of the pushforward operator}
	
	For a Polish space $X$, let $C_\mathsf{b}(X)$ denote the set of all bounded continuous real-valued functions on $X$. By definition, the weak-$\ast$ topology on $\Prob(X)$ is generated by the maps $\Prob(X) \to \R \colon \mu \mapsto \int \xi \,\Diff \mu$, where $\xi \in C_\mathsf{b}(X)$.
	
	\begin{lemma}\label{lemma:push_of_cont}
		Let $X$ and $K$ be Polish spaces and let $f \colon X \to K$ be continuous. Then $f_\ast \colon \Prob(X) \to \Prob(K)$ is also continuous. \hfill \qedsymbol
	\end{lemma}
	
	Now we turn to the continuity properties of the mapping $f \mapsto f_\ast$.
	
	\begin{lemma}\label{lemma:push_mu}
		Let $(K, \dist)$ be a compact metric space and let $(X, \mu)$ be a standard probability space. Then the function
		$
		(\mathfrak{B}(X, K), {\dist}_\mu) \to \Prob(K) \colon f \mapsto f_\ast(\mu)
		$
		is continuous.
	\end{lemma}
	\begin{coolproof}
		Let $f$, $f_0$, $f_1$, \ldots\ $\in \mathfrak{B}(X, K)$ be such that $f_n \to f$ in $(\mathfrak{B}(X, K), \dist_\mu)$.
		To demonstrate that $(f_n)_\ast(\mu) \to f_\ast(\mu)$, let $\xi \in C_\mathsf{b}(X)$; we have to show that
		\begin{equation}\label{eq:xi}
		\int_K \xi \,\Diff(f_n)_\ast(\mu) \,\to\, \int_K \xi \,\Diff f_\ast(\mu).
		\end{equation}
		We may scale $\xi$ if necessary to make it bounded in absolute value by $1$. Take any $\epsilon > 0$. Since $K$ is compact, $\xi$ is uniformly continuous, so we can let $\delta > 0$ be such that $\xi(a) \approx_\epsilon \xi(b)$ whenever $a \approx_\delta b$. For $n \in \N$, let $X_n$ denote the set of all $x \in X$ with $\dist(f_n(x), f(x)) < \delta$. Since $\dist_\mu(f_n, f) \to 0$, we have $\mu(X_n) \to 1$, and hence, for all large enough $n \in \N$,
		\[
		\left|\int_K \xi \,\Diff(f_n)_\ast(\mu) \,-\, \int_K \xi \,\Diff f_\ast(\mu)\right| \,\leq\, \int_X |\xi \circ f_n \,-\, \xi \circ f|\,\Diff\mu \,\leq\, \epsilon \mu(X_n) \,+\, 2(1 - \mu(X_n)) \,\leq\, 2\epsilon.
		\]
		Since $\epsilon$ was chosen arbitrarily, \eqref{eq:xi} follows.
	\end{coolproof}
	
	If $K$ is a compact metric space, then the space $C_\mathsf{b}(K)$, equipped with the uniform norm, is separable. Therefore, there exists a countable set $\set{\xi_n}_{n \in \N}$ of continuous real-valued functions on $K$ bounded in absolute value by $1$ such that $\set{a\xi_n \,:\, a \in \R,\ n\in \N}$ is a dense subset of $C_\mathsf{b}(K)$. With this choice of  $\set{\xi_n}_{n \in \N}$, the topology on $\Prob(K)$ is induced by the metric $\Delta^K$:
	\[
		\Delta^K(\mu, \nu) \defeq \sum_{n = 0}^\infty \frac{\left|\int_K \xi_n \,\Diff \mu \,-\, \int_K \xi_n \,\Diff \nu\right|}{2^{n+1}}.
	\]
	
	\begin{lemma}\label{lemma:push_uni}
		Let $X$ be a Polish space and let $(K, \dist)$ be a compact metric space. Then the map $(\mathfrak{B}(X, K), \dist_{\mathsf{uni}}) \to (\mathfrak{B}(\Prob(X), \Prob(K)), \Delta^K_\mathsf{uni}) \colon f \mapsto f_\ast$ is continuous.
	\end{lemma}
	\begin{coolproof}
		Let $\set{\xi_n}_{n \in \N}$ be the set of functions used to define $\Delta^K$. Take any $N \in \N^+$ and $\epsilon > 0$. Since $K$ is compact, each $\xi_n$ is uniformly continuous, hence we can choose $\delta > 0$ so that $\xi_n(a) \approx_\epsilon \xi_n(b)$ for all $n \leq N$, whenever $a \approx_\delta b$. Let $f$, $g \in \mathfrak{B}(X, K)$ and suppose that $\dist_\mathsf{uni}(f, g) < \delta$. Then, for any $\mu \in \Prob(K)$, we have
		\begin{align*}
			\Delta^K(f_\ast(\mu), g_\ast(\mu)) \,&=\, \sum_{n = 0}^\infty \frac{\left|\int_K \xi_n \,\Diff f_\ast(\mu) \,-\, \int_K \xi_n \,\Diff g_\ast(\mu)\right|}{2^{n+1}} \\
			&\leq\, \sum_{n = 0}^{N} \frac{1}{2^{n+1}}\int_X \left|\xi_n \circ f - \xi_n \circ g\right| \,\Diff \mu \,+\, \frac{1}{2^{N-1}} \,<\, \epsilon + \frac{1}{2^{N-1}}.
		\end{align*}
		Hence, $\Delta^K_\mathsf{uni}(f_\ast, g_\ast) < \epsilon + 2^{-N+1}$. Since $\epsilon$ and $N$ are arbitrary, this completes the proof.
	\end{coolproof}
	
	\subsubsection*{Density of continuous functions}
	
	Recall that a topological space $X$ is \emphd{zero-dimensional} if it has a basis consisting of clopen sets.
	
	\begin{lemma}\label{lemma:zero-dim}
		Let $X$ be a zero-dimensional Polish space and let $(K, \dist)$ be a compact metric space. If $\mu \in \Prob(X)$, then the set of all continuous maps $f \colon X \to K$ is dense in $(\mathfrak{B}(X, K), \dist_\mu)$.
	\end{lemma}
	\begin{coolproof}
		Without loss of generality, assume that the metric $\dist$ is bounded by $1$. Let $f \in \mathfrak{B}(X, K)$ and $\epsilon > 0$. Since $K$ is compact, it contains a finite $\epsilon$-net $Z = \set{z_0, \ldots, z_{n-1}} \subseteq K$. Let $g \colon X \to Z$ be the map that sends each $x \in X$ to the point $z \in Z$ that is closest to $f(x)$ \ep{ties may be broken arbitrarily}. Then $\dist_\mathsf{uni}(f, g) < \epsilon$ by construction. Since the measure $\mu$ is regular \cite[Theorem 17.10]{K_DST} and the space $X$ is zero-dimensional, for each $0 \leq i < n$, there is a clopen set $U_i \subseteq X$ such that $\mu(U_i \symdif g^{-1}(z_i)) < \epsilon/n$. For every $x \in X$, set $h(x) \defeq z_i$ if $x \in V_i \defeq U_i \setminus (U_0 \cup \ldots \cup U_{i-1})$ for some $0 \leq i < n$, and $h(x) \defeq z_0$ if $x \in V \defeq X \setminus (U_0 \cup \ldots \cup U_{n-1})$. Since the sets $V_0$, \ldots, $V_{n-1}$, $V$ are clopen, the map $h$ is continuous. If $h(x) \neq g(x)$, then either $x \in V_i \setminus g^{-1}(z_i)$ for some $0 \leq i < n$, in which case $x \in U_i \setminus g^{-1}(z_i)$; or else, $x \in V \setminus g^{-1}(z_0)$, in which case $x \in g^{-1}(z_i) \setminus U_i$ for some $1 \leq i< n$. Since the metric $\dist$ is bounded by $1$, we conclude that 
		\[
			\dist_\mu(g, h) \,\leq\, \sum_{i = 0}^{n-1} \mu(U_i \setminus g^{-1}(z_i)) \,+\, \sum_{i=1}^{n-1} \mu(g^{-1}(z_i) \setminus U_i) \,\leq\, \sum_{i=0}^{n-1} \mu(U_i \symdif g^{-1}(z_i)) \,<\, \epsilon.
		\]
		Therefore, we have found a continuous function $h \colon X \to K$ with $\dist_\mu(f, h) < 2\epsilon$. As $\epsilon$ is arbitrary, the proof is complete.
	\end{coolproof}
	
	\section{Combinatorial reductions}\label{sec:red}
	
	For $k \in \N^+$, let $u_k$ denote the uniform probability measure on $k$, i.e., let $u_k(\set{i}) \defeq 1/k$ for all $i < k$. Set $\Omega_k \defeq k^\G$ and $\bm{u}_k \defeq u_k^\G$. Recall that the space $(\Omega_k, \bm{u}_k)$ is equipped with the shift action $\sigma_k$.
	
	Given $\phi \in \finf{\G}{k}$ and a partial map $c \colon \G \rightharpoonup k$, we say that $\gamma \in \G$ is an \emphd{occurrence} of $\phi$ in $c$ if $\gamma \cdot c \supseteq \phi$. The set of all {occurrences} of $\phi$ in $c$ is denoted by $\mathcal{O}_\phi(c)$. By definition, if $\gamma \in \mathcal{O}_\phi(c)$, then, in particular, $\dom(\phi) \gamma \subseteq \dom(c)$. Define
	\[
		\Omega_k(\phi) \defeq \set{c \in \Omega_k \,:\, \mathbf{1} \in \mathcal{O}_\phi(c)} = \set{c \in \Omega_k \,:\, \phi \subset c}.
	\]
	Note that $\bm{u}_k(\Omega_k(\phi)) = k^{-|\phi|}$ \ep{where $|\phi|$ is the cardinality of the domain of $\phi$}. The family of sets $\set{\Omega_k(\phi) \,:\, \phi \in \finf{\G}{k}}$ forms a basis for the topology on $\Omega_k$ consisting of clopen sets. From this fact and \cite[Theorem~17.20]{K_DST}, we obtain the following:
	
	\begin{lemma}\label{lemma:basic}
		Let $k \in \N^+$ and $\mu$, $\mu_0$, $\mu_1$, \ldots\ $\in \Prob(\Omega_k)$. Then $\lim_{n \to \infty}\mu_n = \mu$ if and only if, for all $\phi \in \finf{\G}{k}$, we have $\lim_{n \to \infty} \mu_n(\Omega_k(\phi)) \,=\, \mu(\Omega_k(\phi))$. \hfill \qedsymbol
	\end{lemma}
	
	We also consider the space $\tilde{\Omega}_k \defeq (k^\N)^\G$, equipped with the product measure $\tilde{\bm{u}}_k \defeq (u_k^\N)^\G$ and
	the shift action $\sigma_{k^\N}$ of $\G$. To simplify the notation, given $x \in \tilde{\Omega}_k$, $\gamma \in \G$, and $n \in \N$, we write $x(\gamma, n)$ to mean $x(\gamma)(n)$ \ep{however, $x(\gamma)$ still denotes the corresponding element of $k^\N$}. If $k \geq 2$, then, by the {measure isomorphism theorem} \cite[Theorem~17.41]{K_DST}, the standard probability spaces $([0;1], \lambda)$ and $(k^\N, u_k^\N)$ are Borel isomorphic, which allows us to replace $\sigma \colon \G \acts (\Omega, \LG)$ by $\sigma_{k^\N} \colon \G \acts (\tilde{\Omega}_k, \tilde{\bm{u}}_k)$ in the statements of Theorems~\ref{theo:ult_erg}, \ref{theo:main_delta}, and \ref{theo:main_Borel}. This gives us two main advantages. First, the space $\tilde{\Omega}_k$ is zero\=/dimensional; in particular, Lemma~\ref{lemma:zero-dim} applies to it. Second, the structure of $\tilde{\Omega}_k$ will be explicitly used in the proof of Theorem~\ref{theo:ult_erg} presented in \S\ref{subsec:ult_erg}.

	\subsection{Reduction for Theorem~\ref{theo:weak_erg}}
	
	In this subsection we reduce Theorem~\ref{theo:weak_erg} to the following statement:
	
	\begin{theobis}{theo:weak_erg}\label{theo:weak_erg_bis}
		Let $k \in \N^+$ and let $(D_n)_{n \in \N}$ be a sequence of finite subsets of $\G$ with $|D_n|/\log n \to\infty$. Then, for all $S \in \fins{\G}$ and $\phi \colon S \to k$, we have
		\[
			\lim_{n \to \infty} \frac{|D_n \cap \mathcal{O}_\phi(c)|}{|D_n|} \,=\, \frac{1}{k^{|S|}}, \qquad \text{for $\bm{u}_k$-a.e.}\ c \in \Omega_k.
		\]
	\end{theobis}
	
	\begin{lemma}\label{lemma:red1}
		Theorem \ref{theo:weak_erg_bis} implies Theorem~\ref{theo:weak_erg}.
	\end{lemma}
	\begin{coolproof} \stepcounter{ForClaims} \renewcommand{\theForClaims}{\ref{lemma:red1}}
		Assume Theorem~\ref{theo:weak_erg_bis}. Fix a sequence $(D_n)_{n \in \N}$ of nonempty finite subsets of $\G$ such that $|D_n|/\log n \to \infty$. Notice that Theorem~\ref{theo:weak_erg} is equivalent to the following assertion:
		\begin{equation}\label{eq:weak_erg}
			\lim_{n \to \infty} \uni_{x, D_n} = \LG, \qquad \text{for $\LG$-a.e.}\ x \in \Omega.
		\end{equation}
		On the other hand, by Lemma~\ref{lemma:basic}, the conclusion of Theorem~\ref{theo:weak_erg_bis} is equivalent to
		\begin{equation}\label{eq:assump}
			\lim_{n \to \infty} \uni_{c, D_n} \,=\, \bm{u}_k, \qquad \text{for $\bm{u}_k$-a.e.}\ c \in \Omega_k.
		\end{equation}
		
		\begin{claim}\label{claim:factor}
			If $\pi \colon (\Omega, \LG) \to (\Omega_k, \bm{u}_k)$ is a factor map, then
			\[
				\lim_{n \to \infty} \M_{D_n} \pi(x) \,=\, \M_\LG \pi \,=\, \bm{u}_k, \qquad \text{for $\LG$-a.e.}\ x \in \Omega.
			\]
		\end{claim}
		\begin{claimproof}
			From the equivariance of $\pi$, it follows that for all $x \in \Omega$ and $D \in \fins{\G} \setminus \set{\0}$,
			\[
				\M_D \pi(x) \,=\, \pi_\ast(\uni_{x, D}) \,=\, \uni_{\pi(x), D}.
			\]
			Using \eqref{eq:assump} and the fact that, since $\pi$ is a factor map, $\M_\LG(\pi) = \pi_\ast(\LG) = \bm{u}_k$, we conclude that
			\[
				\M_{D_n} \pi(x) \,=\, \uni_{\pi(x), D_n} \,\xrightarrow[n \to \infty]{} \bm{u}_k, \qquad \text{for $\LG$-a.e.}\ x \in \Omega. \qedhere
			\]
		\end{claimproof}

		Define a function $\pr \colon \Omega \to [0;1]$ by $\pr(x) \defeq x(\mathbf{1})$. Notice that $\code{\pr} = \id_\Omega$. For each $k \in \N^+$, let $Q_k \subset [0;1]$ be the set of all fractions of the form $i/k$, $0 \leq i < k$,
			and define $f_k \colon \Omega \to Q_k$ by
			\[
				f_k(x) \defeq \max \set{q \in Q_k \,:\, q \leq \pr(x)}.
			\]
			Let $\pi_k \defeq \code{f_k}$. By definition, $f_k(x) \approx_{1/k} \pr(x)$ for all $x \in \Omega$; in other words, the sequence $(f_k)_{k \in \N}$ converges to $\pr$ uniformly. By Lemmas~\ref{lemma:cont_of_code}\ref{item:uni} and \ref{lemma:push_uni}, this implies that
			\[
				\pi_k \to \id_\Omega \quad \text{and} \quad (\pi_k)_\ast \to \id_{\Prob(\Omega)} \qquad \text{uniformly}.
			\]
			By construction, $(f_k)_\ast(\LG)$ is the uniform probability measure on $Q_k$, and $(\pi_k)_\ast(\LG)$ is the corresponding product measure on $Q_k^\G$. Thus, we may apply Claim~\ref{claim:factor} to $\pi_k$ and conclude that
			\[
				\lim_{n \to \infty} \M_{D_n} \pi_k(x) \,=\, \M_\LG \pi_k, \qquad \text{for $\LG$-a.e.}\ x \in \Omega.
			\]
			We can put all of these facts together as follows:
			
			\begin{figure}[H]
			\begin{tikzpicture}
				\node (a) at (0,0) {$\M_{D_n} \pi_k (x)$};
				\node (b) at (0,-2.3) {$\M_\LG \pi_k$};
				\node (c) at (5,0) {$\uni_{x, D_n}$};
				\node (d) at (5,-2.3) {$\LG$};
				
				\draw[->] (a) -- node[midway, anchor=south, rotate=-90] {$n \to \infty$} (b);
				\draw[->] (b) -- node[midway, anchor=south] {$k \to \infty$} (d);
				\draw[->] (a) -- node[midway, anchor=south] {$k \to \infty$} node[midway, anchor=north] {uniformly in $n$} (c);
			\end{tikzpicture}
			\end{figure}
			
			\noindent It is clear from the above diagram that $\uni_{x, D_n}$ converges to $\LG$ as $n \to \infty$, proving \eqref{eq:weak_erg}.
	\end{coolproof}

	\subsection{Reductions for Theorems~\ref{theo:ult_erg}, \ref{theo:main_delta}, and \ref{theo:main_Borel}}
	
	Theorem~\ref{theo:ult_erg} reduces to the following statement:
	
	\begin{theobis}{theo:ult_erg}\label{theo:ult_erg_bis}
		For all $k \in \N^+$, $S \in \fins{\G}$, and $\epsilon > 0$, there is $C > 0$ with the following property:
		
		\smallskip
		
		Let $(D_n)_{n \in \N}$ be a sequence of finite subsets of $\G$ with $|D_n| \geq C \log(n+2)$ for all $n \in \N$. Then there exists a Borel map $g \colon \tilde{\Omega}_k \to k$ such that
		\[
		\tilde{\bm{u}}_k(\set{x \in \tilde{\Omega}_k \,:\, g(x) \neq x(\mathbf{1}, 0)}) \,\leq\, \epsilon,
		\]
		and, for all $\phi \colon S \to k$, we have
		\[
		\frac{|D_n \cap \mathcal{O}_\phi(\code{g}(x))|}{|D_n|}\,\approx_\epsilon\, \frac{1}{k^{|S|}}, \qquad \text{for all $n \in \N$ and for $\tilde{\bm{u}}_k$-a.e.}\ x \in \tilde{\Omega}_k.
		\]
	\end{theobis}
	
	\begin{lemma}\label{lemma:red2}
		Theorem~\ref{theo:ult_erg_bis} implies Theorem~\ref{theo:ult_erg}.
	\end{lemma}
	\begin{coolproof}\stepcounter{ForClaims} \renewcommand{\theForClaims}{\ref{lemma:red2}}
		Assume Theorem~\ref{theo:ult_erg_bis}. Taking advantage of the measure isomorphism theorem, we will prove the statement of Theorem~\ref{theo:ult_erg} with $(\tilde{\Omega}_2, \tilde{\bm{u}}_2)$ in place of $(\Omega, \LG)$. We equip the Cantor space $2^\N$ with the metric $\mathfrak{m}$ given by
		\[
			\mathfrak{m}(a, b) \defeq \sum_{n = 0}^\infty \frac{\mathbbm{1}_{a(n) \neq b(n)}}{2^{n+1}}.
		\]
		Define $\pr \colon \tilde{\Omega}_2 \to 2^\N$ by $\pr(x) \defeq x(\mathbf{1})$. Note that $\code{\pr} = \id_{\tilde{\Omega}_2}$ and $\M_{\tilde{\bm{u}}_2} \code{\pr} = \tilde{\bm{u}}_2$.
		
		\begin{claim}
			It suffices to prove Theorem~\ref{theo:ult_erg} with $K = 2^\N$ and $f = \pr$.
		\end{claim}
		\begin{claimproof}
			Let $(K, \dist)$ be a compact metric space. Without loss of generality, assume that the metric $\dist$ is bounded by $1$. Fix a Borel function $f \colon \tilde{\Omega}_2 \to K$, $\epsilon > 0$, and an open neighborhood $U$ of $\M_{\tilde{\bm{u}}_2}\code{f}$. The space $\tilde{\Omega}_2$ is zero-dimensional, so Lemmas \ref{lemma:zero-dim}, \ref{lemma:cont_of_code}\ref{item:mu}, and \ref{lemma:push_mu} allow us to assume that $f$ is continuous \ep{after replacing $\epsilon$ by, say, $\epsilon/2$}. Since $\tilde{\Omega}_2$ is compact, $f$ is uniformly continuous, so we can pick $\delta \in (0; \epsilon/2)$ such that $f(x) \approx_{\epsilon/2} f(y)$ whenever $x \approx_\delta y$. By Lemma~\ref{lemma:push_of_cont}, the set $U' \defeq (\code{f})_\ast^{-1}(U)$ is an open neighborhood of $\tilde{\bm{u}}_2$.
			
			Let $\mathsf{q} \colon \tilde{\Omega}_2 \to 2^\N$ be a Borel map and consider the function $g \defeq f \circ \mathsf{q}$. Note that if $\mathfrak{m}_{\tilde{\bm{u}}_2}(\pr, \mathsf{q}) \leq \delta^2$, then $\dist_{\tilde{\bm{u}}_2}(f, g) \leq \epsilon$. Indeed, if $\mathfrak{m}_{\tilde{\bm{u}}_2}(\pr, \mathsf{q}) \leq \delta^2$, then, by Markov's inequality,
			\[
				\tilde{\bm{u}}_2(\set{x \in \tilde{\Omega}_2 \,:\, \pr(x) \not \approx_\delta \mathsf{q}(x)}) \,\leq\, \delta,
			\]
			and, by the choice of $\delta$ and since $\dist$ is bounded by $1$, we have $\dist_{\tilde{\bm{u}}_2}(f, g) \leq \epsilon/2 + \delta < \epsilon$. Additionally, if $D\in \fins{\G}\setminus \set{\0}$ and $x \in \tilde{\Omega}_2$ satisfy $\M_D \code{\mathsf{q}}(x) \in U'$,
			then
			\[
				\M_{D} \code{g}(x) \,=\, (\code{f})_\ast (\M_{D} \code{\mathsf{q}}(x)) \,\in \, U.
			\]
			Therefore, if Theorem~\ref{theo:ult_erg} holds for $\pr$, $\delta^2$, and $U'$, then it also holds for $f$, $\epsilon$, and $U$, as desired.
		\end{claimproof}
		
		The remainder of the argument is similar to the last part of the proof of Lemma~\ref{lemma:red1}. For each $n \in \N^+$, let $Q_n$ be the set of all $a \in 2^\N$ such that $a(i) = 0$ for all $i \geq n$, and define $\pr_n \colon \tilde{\Omega}_2 \to Q_n$ by
		\[
			\pr_n(x)(i) \defeq \begin{cases}
				x(\mathbf{1}, i) &\text{if } i < n;\\
				0 &\text{if } i \geq n.
			\end{cases}
		\]
		Then $\pr_n \to \pr$ uniformly, so to prove Theorem~\ref{theo:ult_erg} for $\pr$, it is enough to prove it for each $\pr_n$. Due to Lemma~\ref{lemma:basic}, Theorem~\ref{theo:ult_erg} for $\pr_1$ is equivalent to Theorem~\ref{theo:ult_erg_bis} applied with $k = 2$. For larger $n$, consider the mapping $\theta_n \colon 2^\N \to (2^n)^\N$ given by
		\[
			\theta(a)(i) \defeq (a(i  n),\ a(i  n + 1),\ \ldots,\ a(i  n + n - 1)) \qquad \text{for all } a \in 2^\N \text{ and } i \in \N,
		\]
		where we identify the natural numbers less than $2^n$ with the $n$-tuples of zeros and ones. This mapping induces an equivariant isomorphism between $(\tilde{\Omega}_2, \tilde{\bm{u}}_2)$ and $(\tilde{\Omega}_{2^n}, \tilde{\bm{u}}_{2^n})$ and shows that Theorem~\ref{theo:ult_erg} for $\pr_n$ is equivalent to Theorem~\ref{theo:ult_erg_bis} applied with $k = 2^n$.
	\end{coolproof}
	
	Similarly, Theorems \ref{theo:main_delta} and \ref{theo:main_Borel} reduce to the following statements:
	
	\begin{theobis}{theo:main_delta}\label{theo:main_delta_bis}
		For all $k \in \N^+$, $S \in \fins{\G}$, and $\epsilon > 0$, there is $C > 0$ with the following property:
		
		\smallskip
		
		Let $D$ be a finite subset of $\G$ with $|D| \geq C$ and let $\alpha \colon \G \acts X$ be an $(S, D)$-free Borel action of $\G$. Then, for any $\mu \in \Prob(X)$ and $\delta > 0$, there is a Borel map $g \colon X \to k$ such that, for all $\phi\colon S \to k$,
		\[
			\mu\left(\left\{x \in X \,:\, \frac{|D \cap \mathcal{O}_\phi(\code{g}(x))|}{|D|} \approx_\epsilon \frac{1}{k^{|S|}}\right\}\right) \,\geq\, 1 - \delta.
		\]
	\end{theobis}
	
	\begin{theobis}{theo:main_Borel}\label{theo:main_Borel_bis}
		For all $k \in \N^+$, $S \in \fins{\G}$, and $\epsilon > 0$, there is $C > 0$ with the following property:
		
		\smallskip
		
		Let $D$ be a finite subset of $\G$ with $|D| \geq C$ and let $\alpha \colon \G \acts X$ be a uniformly subexponential $(S, D)$-free Borel action of $\G$. Then there is a Borel map $g \colon X \to k$ such that, for all $\phi \colon S \to k$,
		\[
		\frac{|D \cap \mathcal{O}_\phi(\code{g}(x))|}{|D|} \,\approx_\epsilon\, \frac{1}{k^{|S|}} \qquad \text{for all } x \in X.
		\]
	\end{theobis}
	
	The proof of the following lemma is essentially the same as of Lemma~\ref{lemma:red2}, and we omit it.
	
	\begin{lemma}
		Theorem~\ref{theo:main_delta_bis} implies Theorem~\ref{theo:main_delta}, while Theorem~\ref{theo:main_Borel_bis} implies Theorem~\ref{theo:main_Borel}.	\hfill \qedsymbol
	\end{lemma}
	
	\section{Using concentration of measure}\label{sec:concentration}
	
	\subsection{The main probabilistic bound}
	
	The following inequality is the main probabilistic input for our arguments:
	
	\begin{lemma}\label{lemma:concentrated}
		Let $k \in \N^+$, $S \in \fins{\G}$, and $\epsilon > 0$. Let $D$ be a nonempty finite subset of $\G$ and let $\alpha \colon \G \acts X$ be an $(S, D)$-free action of $\G$. Take any $x \in X$ and pick a function $c \colon (SD \cdot x) \to k$ uniformly at random. Then, for all $\phi \colon S \to k$,
		\[
			\P \left[\frac{|D \cap \mathcal{O}_\phi(\code{c}(x))|}{|D|} \not\approx_\epsilon \frac{1}{k^{|S|}}\right] \,\leq\, 2\exp\left(-\epsilon^2 \frac{|D|}{2 |S|^3}\right).
		\]
	\end{lemma}
	\begin{coolproof}
		We shall apply the following concentration of measure result, which is a consequence of Azuma's inequality for Doob martingales \cite[\S{}7.4]{AS}:
		
		\begin{theo}[{\textls{Simple Concentration Bound}; see \cite[79]{MR}}]\label{theo:SCB1}
			Let $\xi$ be a random variable determined by $s$ independent trials such that changing the outcome of any one trial can affect $\xi$ at most by~$b$. Then 
			\[
			\P\left[\xi \not \approx_t \mathbb{E} \xi\right] \,\leq\, 2\exp\left(-\frac{t^2}{2b^2s}\right).
			\]
		\end{theo}
		
		\noindent Let $k$, $S$, $\epsilon$, $D$, $\alpha$, and $c$ be as in the statement of Lemma~\ref{lemma:concentrated}. Since the action $\alpha$ is $S$-free, for all $\phi \colon S \to k$, we have
		\[
		\mathbb{E}\left[\left|D \cap \mathcal{O}_\phi(\code{c}(x))\right|\right] \,=\, \sum_{\delta \in D} \P\left[\delta \in \mathcal{O}_\phi(\code{c}(x))\right] \,=\, \frac{|D|}{k^{|S|}}.
		\]
		Consider any $y \in SD \cdot x$ and let $c_1$, $c_2 \colon (SD\cdot x) \to k$ be two maps that agree on $(SD \cdot x) \setminus \set{y}$. Let $\phi \colon S \to k$ and suppose that some $\delta \in D$ belongs to $\mathcal{O}_\phi(\code{c_1}(x)) \symdif \mathcal{O}_\phi(\code{c_2}(x))$. Then $y \in S \cdot (\delta \cdot x)$, i.e., $\delta \cdot x \in S^{-1} \cdot y$. Since $\alpha$ is $D$-free, there are at most $|S^{-1} \cdot y| = |S|$ possible values for $\delta$. Thus, we may apply the \hyperref[theo:SCB1]{Simple Concentration Bound} with parameters
		\[
		s \defeq |SD \cdot x| \leq |S||D|, \qquad b \defeq |S|, \qquad \text{and} \qquad t \defeq \epsilon|D|,
		\]
		to obtain
		\[
		\P \left[|D \cap \mathcal{O}_\phi(\code{c}(x))| \not\approx_{\epsilon |D|} \frac{|D|}{k^{|S|}}\right] \,\leq\, 2\exp\left(-\epsilon^2 \frac{|D|}{2 |S|^3}\right),
		\]
		as desired.
	\end{coolproof}
	
	\subsection{Proof of Theorem~\ref{theo:weak_erg}}\label{subsec:weak_erg}
	
	We are now ready to prove Theorem~\ref{theo:weak_erg} \ep{or rather Theorem~\ref{theo:weak_erg_bis}}.
		Let $k \in \N^+$ and let $(D_n)_{n \in \N}$ be a sequence of finite subsets of $\G$ such that $|D_n|/\log n \to\infty$. Take any $S \in \fins{\G}$, $\phi \colon S \to k$, and $\epsilon > 0$. We will show that for $\bm{u}_k$-a.e.\ $c \in \Omega_k$ and for all sufficiently large $n \in \N$,
	\begin{equation}\label{eq:approx_epsilon}
		\frac{|D_n \cap \mathcal{O}_\phi(c)|}{|D_n|} \,\approx_\epsilon\, \frac{1}{k^{|S|}},
	\end{equation}
	which will imply the conclusion of Theorem~\ref{theo:weak_erg_bis}. For each $n \in \N$, let $X_n$ denote the set of all $c \in \Omega_k$ for which \eqref{eq:approx_epsilon} fails. By Lemma~\ref{lemma:concentrated}, we have
	\[
		\sum_{n \in \N} \bm{u}_k(X_n) \,\leq\, \sum_{n \in \N} 2\exp\left(-\epsilon^2 \frac{|D_n|}{2 |S|^3}\right) \,<\, \infty,
	\]
	since $\epsilon^2|D_n|/(2|S|^3) > 2 \log n$ for all sufficiently large $n$. An application of the Borel--Cantelli lemma completes the proof.
	
	\section{The Lov\'asz Local Lemma and its measurable versions}\label{sec:LLL}
	
	\subsection{The classical LLL}\label{subsec:SLLL}
	
	The reader is referred to~\cite[Chapter~5]{AS} and \cite{MR} for background on the~LLL and its applications in combinatorics. The presentation below follows, with slight modifications,~\cite[Section~1.2]{MLLL}.
	
	Let $X$ be a set and let $k \in \N^+$. A \emphd{bad \ep{$k$-}event} over $X$ is a nonempty subset $B \subseteq \finf{X}{k}$ such that for all $\phi$, $\phi' \in B$, $\dom(\phi) = \dom(\phi')$. If a bad event $B$ is nonempty, then its \emphd{domain} is the set $\dom(B) \coloneqq \dom(\phi)$ for any \ep{hence all} $\phi \in B$; the domain of the empty bad event is, by definition, the empty set. The \emphd{probability} of a bad $k$-event $B$ with domain $F$ is defined to be
	\[
		\P[B] \coloneqq {|B| \over k^{|F|}}.
	\]
	We say that a map $f \colon X \to k$ \emphd{avoids} a bad $k$-event $B$ if there is no $\phi \in B$ such that $\phi \subseteq f$. Note that if $X$ is finite and $f \colon X \to k$ is chosen uniformly at random, then $\P[B]$ is the probability that $f$ does not avoid $B$.
	
	A~\emphd{\ep{$k$\=/}instance \ep{of the LLL}} over a set $X$ is an arbitrary set $\B$ of bad $k$-events. A~\emphd{solution} to a $k$-instance~$\B$ is a function $f \colon X \to k$ that avoids all $B \in \B$. For an instance $\B$ and a bad event $B \in \B$, the \emphd{neighborhood} of $B$ in $\B$ is the set
	\[
		\Nbhd_\B(B) \defeq \set{B' \in \B \setminus \set{B} \,:\, \dom(B') \cap \dom(B) \neq \0}.
	\]
	The \emphd{degree} of $B$ in $\B$ is defined to be $\deg_\B(B) \defeq |\Nbhd_\B(B)|$. Let
	\[
		p(\B) \defeq \sup_{B \in \B} \P[B] \qquad \text{and} \qquad d(\B) \defeq \sup_{B \in \B} \deg_\B(B).
	\]
	An instance $\B$ is \emphd{correct for the Symmetric LLL} (the \emphd{SLLL} for short) if
	\[
		e\cdot p(\B) \cdot (d(\B) + 1) < 1,
	\]
	where $e = 2.71\ldots${} denotes the base of the natural logarithm. Note that if $\B$ is correct for the SLLL, then, in particular, $\deg_\B(B) < \infty$ for all $B \in \B$ \ep{instances $\B$ with this property are called \emphd{locally finite} in \cite{MLLL}}.
	
	\begin{theo}[{Erd\H os--Lov\'asz~\cite{EL}; \textls{Symmetric Lov\'asz Local Lemma}}]\label{theo:SLLL}
		Let $k \in \N^+$ and let $\B$ be a $k$-instance of the LLL over a set $X$. If $\B$ is correct for the SLLL, then $\B$ has a solution.
	\end{theo}
	
	The Symmetric LLL was introduced by Erd\H os and Lov\'asz \ep{with $4$ in place of $e$} in their seminal paper~\cite{EL}; the constant was subsequently improved by Lov\'asz (the sharpened version first appeared in~\cite{Spencer}). Theorem~\ref{theo:SLLL} is a special case of the~SLLL in the so-called \emph{variable framework} (the~name is due to Kolipaka and Szegedy~\cite{KolipakaSzegedy}), which encompasses most typical applications. For the full statement of the~SLLL, see~\cite[Corollary~5.1.2]{AS}. Deducing Theorem~\ref{theo:SLLL} from \cite[Corollary~5.1.2]{AS} is routine when $X$ is finite (see, e.g., \cite[41]{MR}); the case of infinite~$X$ then follows by compactness. A more general version of Theorem~\ref{theo:SLLL} for infinite $X$, with $k$ replaced by an arbitrary standard probability space, was proved by Kun~\cite[Lemma~13]{Kun}.
	
	Theorem~\ref{theo:SLLL} can be extended to instances $\B$ with $d(\B) = \infty$, provided that the probability $\P[B]$ of a bad event $B \in \B$ decays sufficiently quickly as $|\dom(B)|$ increases. An instance $\B$ is \emphd{correct for the General LLL} (the~\emphd{GLLL} for short) if $\Nbhd_\B(B)$ is countable for every $B \in \B$, and there is a function $\omega \colon \B \to [0;1)$, called a \emphd{witness} to the correctness of $\B$, such that for all $B \in \B$,
	\[
	\P[B] \leq \omega(B) \prod_{B' \in \Nbhd_\B(B)} (1 - \omega(B')).
	\]
	
	\begin{theo}[{\textls{General Lov\'asz Local Lemma}; \cite[Lemma~5.1.1]{AS}}]\label{theo:GLLL}
		Let $k \in \N^+$ and let $\B$ be a $k$-instance of the~LLL over a set $X$. If $\B$ is correct for the~GLLL, then $\B$ has a solution.
	\end{theo}
	
	A standard calculation (see~\cite[proof of Corollary~5.1.2]{AS}) shows that if an instance $\B$ is correct for the~SLLL, then it is also correct for the~GLLL (hence the name ``General LLL'').
	
	\subsection{Measurable versions of the LLL}\label{subsec:MLLL}
	
	Let $X$ be a standard Borel space and let $k \in \N^+$. Then the set of all bad $k$-events is also naturally equipped with the structure of a standard Borel space (indeed, each bad event is a finite set, so the set of all bad $k$-events is a Borel subset of the space $\fins{\finf{X}{k}}$). Thus, it makes sense to talk about \emphd{Borel} instances of the~LLL, i.e., Borel sets of bad events.
	
	Given a Borel $k$-instance $\B$ over $X$ that is correct for the~SLLL, it is natural to wonder if it has a Borel solution $f \colon X \to k$. Although the answer is negative in general (see~\cite[Theorem~1.6]{CJMST-D}), Cs\'oka, Grabowski, M\'ath\'e, Pikhurko, and Tyros~\cite{CGMPT} answered the question in the affirmative for \emph{uniformly subexponential} instances. Given an instance $\B$ over a set $X$, an element $x \in X$, and an integer $n \in \N$, let $R^n_\B(x)$ denote the set of all $y \in X$ such that either $y = x$, or there exists a sequence $B_1$, \ldots, $B_m \in \B$ with $m \leq n$ satisfying
	\[
	x \in \dom(B_1), \qquad \dom(B_i) \cap \dom(B_{i+1}) \neq \0 \text{ for all } 1 \leq i < m,\qquad\text{ and }\qquad y \in \dom(B_m).
	\]
	The instance $\B$ is \emphd{uniformly subexponential} if for every $\epsilon > 0$, there exists $n_0 \in \N$ such that for all $n \geq n_0$ and for all $x \in X$, $|R^n_\B(x)| < (1+\epsilon)^n$.
	
	\begin{theo}[{Cs\'oka--Grabowski--M\'ath\'e--Pikhurko--Tyros~\cite[Theorem~1.3]{CGMPT}, \textls{Borel SLLL for uniformy subexponential instances}}]\label{theo:LLL_Borel}
		Let $k \in \N^+$ and let $\B$ be a Borel $k$-instance of the~LLL over a standard Borel space $X$. If $\B$ is correct for the~SLLL and uniformly subexponential, then $\B$ has a Borel solution $f \colon X \to k$.
	\end{theo}
	
	For a $k$-instance $\B$ over a set $X$ and a map $f \colon X \to k$, we define the \emphd{defect} $\Def(f; \B)$ of $f$ with respect to $\B$ by
	\begin{equation}\label{eq:defect1}
	\Def(f,\B) \defeq \set{x \in X \,:\, x \in \dom(\phi) \text{ for some } \phi \in B \in \B \text{ with } \phi \subseteq f}.
	\end{equation}
	Evidently, $f$ is a solution to $\B$ if and only if $\Def(f, \B) = \0$. Thus, in the absence of a Borel solution to $\B$, it is natural to seek a Borel map $f \colon X \to k$ whose defect is ``small'' in some sense. The next result was proved by the current author in~\cite{MLLL}:
	
	\begin{theo}[{\cite[Theorem~5.1]{MLLL}, \textls{approximate SLLL}}]\label{theo:LLL_approx}
		Let $k \in \N^+$ and let $\B$ be a Borel $k$-instance of the~LLL over a standard Borel space $X$. If $\B$ is correct for the~SLLL, then for any $\mu \in \Prob(X)$ and $\delta > 0$, there is a Borel function $f \colon X \to k$ with $\mu(\Def(f, \B)) \leq \delta$.
	\end{theo}
	
	It is an open question whether the conclusion of Theorem~\ref{theo:LLL_approx} holds with $\delta = 0$; see Problem~\ref{prob:meas_LLL} in Section~\ref{sec:remks}. Also, Theorem~\ref{theo:LLL_approx} fails for instances that are correct for the GLLL instead of the SLLL \ep{see \cite[Theorem~7.1]{MLLL} and Remark~\ref{remk:amenable} below}. However, when the underlying structure is in a certain sense induced by the shift action $\sigma$, even instances that are only correct for the~GLLL can be solved with a null defect---see Theorem~\ref{theo:LLL_meas} in the next subsection.
	
	\subsection{Using the LLL over group actions}
	
	Now we describe a convenient set-up for applying the LLL to problems in ergodic theory.
	
	Let $\alpha \colon \G \acts X$ be an action of $\G$ and let $\Phi \subseteq \finf{\G}{k}$ be a bad $k$-event over~$\G$ with domain $F \in \fins{\G}$. For each $x \in X$, define a bad $k$-event $B_x(\Phi, \alpha)$ over $X$ via
	\[
		B_x(\Phi, \alpha) \defeq \set{\phi \colon (F \cdot x) \to k \,:\, \rest{\code{\phi}(x)}{F} \in \Phi}.
	\]
	Note that if $B_x(\Phi, \alpha) \neq \0$, then $\dom(B_x(\Phi, \alpha)) = F \cdot x$. \ep{If $\alpha$ is not $F$-free, then $B_x(\Phi, \alpha)$ may be empty even if $\Phi$ is not.} By construction, a function $f \colon X \to k$ avoids $B_x(\Phi, \alpha)$ precisely when $\code{f}(x)$ avoids $\Phi$. Define an instance $\B(\Phi, \alpha)$ of the LLL over $X$ as follows:
	\[
	\B(\Phi, \alpha) \coloneqq \set{B_x(\Phi, \alpha) \,:\, x \in X}.
	\]
	Clearly, if $X$ is a standard Borel space and $\alpha \colon \G \acts X$ is a Borel action, then the instance $\B(\Phi, \alpha)$ is Borel. A function $f \colon X \to k$ is a solution to $\B(\Phi, \alpha)$ if and only if $\code{f}(x)$ avoids $\Phi$ for all $x \in X$. Hence, it is somewhat more convenient to define the \emphd{defect} of a map $f \colon X \to k$ as the set of all $x \in X$ such that $\code{f}(x)$ does \emph{not} avoid $\Phi$:
	\[
		\Def(f, \Phi, \alpha) \defeq \set{x \in X \,:\, \rest{\code{f}(x)}{F} \in \Phi}.
	\]
	There is a straightforward relationship between this definition and the one in \eqref{eq:defect1}; namely,
	\begin{equation}\label{eq:two_defects}
	\Def(f, \B(\Phi, \alpha)) = F \cdot \Def(f, \Phi, \alpha).
	\end{equation}
	Using the above notation, we can formulate the following corollaries of Theorems~\ref{theo:LLL_Borel} and \ref{theo:LLL_approx}:

	\begin{corl}[to Theorem~\ref{theo:LLL_Borel}]\label{corl:LLL_Borel}
		Let $\alpha \colon \G \acts X$ be a uniformly subexponential Borel action of $\G$ and let $k \in \N^+$. Let $\Phi$ be a bad $k$-event over $\G$ and suppose that the instance $\B(\Phi, \alpha)$ is correct for the SLLL. Then $\B(\Phi, \alpha)$ has a Borel solution $f \colon X \to k$.
	\end{corl}
	
	\begin{corl}[to Theorem~\ref{theo:LLL_approx}]\label{corl:LLL_approx}
		Let $\alpha \colon \G \acts X$ be a Borel action of $\G$ and let $k \in \N^+$. Let $\Phi$ be a bad $k$-event over $\G$ and suppose that the instance $\B(\Phi, \alpha)$ is correct for the SLLL. Then, for any $\mu \in \Prob(X)$ and $\delta > 0$, there is a Borel function $f \colon X \to k$ with $\mu(\Def(f,\Phi, \alpha)) < \delta$.
	\end{corl}
	
	\begin{remk}
		In the statement of Corollary~\ref{corl:LLL_approx}, the measure $\mu$ is not assumed to be $\alpha$-invariant. Because of that, to derive Corollary~\ref{corl:LLL_approx}, one has to apply Theorem~\ref{theo:LLL_approx} not to $\mu$ itself, but to the measure obtained by shifting $\mu$ by one of the elements of $\dom(\Phi)$, and then use \eqref{eq:two_defects}.
	\end{remk}
	
	More generally, let $(\Phi_n)_{n \in \N}$ be a sequence of bad $k$-events over $\G$. For an action $\alpha \colon \G \acts X$ and a map $f \colon X \to k$, define
	\[
	\B((\Phi_n)_{n \in \N}, \alpha) \coloneqq \bigcup_{n = 0}^\infty \B(\Phi_n, \alpha) \qquad \text{and} \qquad \Def(f, (\Phi_n)_{n \in \N}, \alpha) \defeq \bigcup_{n=0}^\infty \Def(f, \Phi_n, \alpha).
	\]
	When $\alpha = \sigma$, we have the following strengthening of Corollary~\ref{corl:LLL_approx}:
	
	\begin{theo}[{\cite[Corollary 6.7]{MLLL}, \textls{measurable GLLL over the shift}}]\label{theo:LLL_meas}
		Let $k \in \N^+$ and let $(\Phi_n)_{n \in \N}$ be a sequence of bad $k$-events over $\G$. Suppose that the instance $\B((\Phi_n)_{n \in \N}, \Free(\sigma))$ is correct for the GLLL. Then there is a Borel function $f \colon \Omega \to k$ with $\LG(\Def(f, (\Phi_n)_{n \in \N}, \sigma)) = 0$.
	\end{theo}
	\begin{remk}\label{remk:amenable}
		Theorem \ref{theo:LLL_meas} can fail for actions other than $\sigma$: According to \cite[Theorem~7.1]{MLLL}, if $\G$ is amenable, then the analog of Theorem~\ref{theo:LLL_meas} holds for a free ergodic \pmp action $\alpha \colon \G \acts (X, \mu)$ if and only if there is a factor map $\pi \colon (X, \mu) \to (\Omega, \LG)$.
	\end{remk}
	
	Theorem~\ref{theo:LLL_meas} is a special case of \cite[Theorem~6.6]{MLLL}, whose full statement is rather technical and will not be needed here. Roughly speaking, \cite[Theorem~6.6]{MLLL} asserts that any combinatorial argument proceeding via a series of iterative applications of the~GLLL can be performed in a measurable fashion over the shift action $\sigma \colon \G \acts (\Omega, \LG)$.

	\section{Proofs of Theorems \ref{theo:ult_erg}, \ref{theo:main_delta}, and \ref{theo:main_Borel}}\label{sec:proofs}
	
	\subsection{Proofs of Theorems \ref{theo:main_delta} and \ref{theo:main_Borel}}\label{subsec:pointwise}
	
	We first establish Theorems~\ref{theo:main_delta} and \ref{theo:main_Borel}, as their proofs are somewhat more straightforward than that of Theorem~\ref{theo:ult_erg} \ep{for instance, they only use the Symmetric LLL rather than the more technical General LLL}.
	
	Let $k \in \N^+$, $S \in \fins{\G}$, and $\epsilon > 0$. For a nonempty finite subset $D \subset \G$, let $\Phi(k, S, \epsilon, D)$ denote the bad $k$-event over $\G$ with domain $SD$ consisting of all maps $c \colon SD \to k$ such that
	\[
	\frac{|D \cap \mathcal{O}_\phi(c)|}{|D|} \,\not\approx_\epsilon\, \frac{1}{k^{|S|}} \qquad \text{for some } \phi \colon S \to k.
	\]
	By definition, if $\alpha \colon \G \acts X$ is a Borel action of $\G$ and $g \colon X \to k$ is a Borel map, then we have
	\begin{equation}\label{eq:defect2}
		x \in \Def(g, \Phi(k, S, \epsilon, D), \alpha) \quad\Longleftrightarrow\quad \frac{|D \cap \mathcal{O}_\phi(\code{g}(x))|}{|D|} \,\not\approx_\epsilon\, \frac{1}{k^{|S|}} \ \text{for some}\  \phi \colon S \to k.
	\end{equation}
	
	\begin{lemma}\label{lemma:correct}
		Let $k \in \N^+$, $S \in \fins{\G}$, and $\epsilon > 0$. There exists $C > 0$ such that for all $D \in \fins{\G}$ with $|D| > C$ and for every $(S, D)$-free action $\alpha \colon \G \acts X$, the instance $\B(\Phi(k, S, \epsilon, D), \alpha)$ is correct for the SLLL.
	\end{lemma}
	\begin{coolproof}
		Let $D \in \fins{\G} \setminus \set{\0}$ and let $\alpha \colon \G \acts X$ be $(S, D)$-free. Set
		\[
		\Phi \defeq \Phi(k, S, \epsilon, D), \qquad \B \defeq \B(\Phi, \alpha),  \qquad \text{and} \qquad B_x \defeq B_x(\Phi, \alpha) \text{ for all } x \in X.
		\]
		Due to Lemma~\ref{lemma:concentrated}, we have
		\[
			p(\B) \,\leq\, 2k^{|S|} \exp\left(-\epsilon^2 \frac{|D|}{2 |S|^3}\right).
		\]
		To upper bound $d(\B)$, note that for each $x \in X$,
		\[
		\Nbhd_\B(B_x) = \set{B_y \in \B \setminus \set{B_x} \,:\, (SD \cdot y) \cap (SD \cdot x) \neq \0}.
		\]
		Since $(SD \cdot y) \cap (SD \cdot x) \neq \0$ if and only if $y \in (SD)^{-1}SD \cdot x$, we obtain
		\[
		\deg_\B(B_x) \leq |(SD)^{-1}SD| - 1 \leq |S|^2|D|^2 - 1.
		\]
		(We subtracted $1$ since $y$ cannot be equal to $x$.) Hence, $d(\B) \leq |S|^2|D|^2 - 1$, and $\B$ is correct for the SLLL as long as
		\[
			e \cdot 2k^{|S|} \exp\left(-\epsilon^2 \frac{|D|}{2 |S|^3}\right) \cdot |S|^2|D|^2 \,<\, 1,
		\]
		which holds whenever $|D|$ is sufficiently large.
	\end{coolproof}
	
	Theorems~\ref{theo:main_delta_bis} and \ref{theo:main_Borel_bis} now follow immediately by combining \eqref{eq:defect2} and Lemma~\ref{lemma:correct} with Corollaries \ref{corl:LLL_approx} and \ref{corl:LLL_Borel} respectively.

	\subsection{Proof of Theorem~\ref{theo:ult_erg}}\label{subsec:ult_erg}
	
	For the purposes of proving Theorem~\ref{theo:ult_erg}, the role of Lemma~\ref{lemma:correct} is played by the following fact:
	
	\begin{lemma}\label{lemma:correct2}
		Let $k \in \N^+$, $S \in \fins{\G}$, and $\epsilon > 0$. There exists $C > 0$ with the following property:
		
		\smallskip
		
		Let $(D_n)_{n \in \N}$ be a sequence of finite subsets of $\G$ with $|D_n|\geq C \log(n+2)$ for all $n \in \N$ and let $\alpha \colon \G \acts X$ be a free action of $\G$. Set
		\[
			\Phi_n \defeq \Phi(k, S, \epsilon, D_n) \text{ for all } n \in \N,
		\]
		\[
			\B \defeq \B((\Phi_n)_{n \in \N}, \alpha), \qquad \text{and} \qquad B_{n,x} \defeq B_x(\Phi_n, \alpha) \text{ for all } n \in \N \text{ and } x \in X.
		\]
		Then the instance $\B$ is correct for the GLLL. Moreover, there is a function $\omega \colon \N \to [0;1)$ such that
		\begin{equation}\label{eq:small_sum}
		\sum_{n = 0}^\infty |SD_n| \cdot \frac{\omega(n)}{1-\omega(n)} \,<\, \epsilon,
		\end{equation}
		and the mapping $\tilde{\omega} \colon \B \to [0;1) \colon B_{n,x} \mapsto \omega(n)$ is a witness to the correctness of $\B$.
	\end{lemma}
	\begin{coolproof}
		Fix any $0 < a < \epsilon^2/(2|S|^3)$. We claim that if $C$ is large enough, then the function \[\omega(n) \defeq \exp(-a|D_n|)\] has the desired properties. To begin with, we are going to assume that $C$ is so large that
		\begin{equation*}
			\exp(-a \cdot C\log(2)) \,<\, 1/2,
		\end{equation*}
		and that the function $\xi \mapsto \xi\exp(-a\xi)$ is decreasing for all $\xi \geq C\log 2$. For any  such $C$, we have
		\begin{align*}
			\sum_{n =0}^\infty |SD_n| \cdot \frac{\omega(n)}{1 - \omega(n)} \,&\leq\, |S|\sum_{n=0}^\infty |D_n| \cdot \frac{\exp(-a|D_n|)}{1 - \exp(-a|D_n|)} \\
			&\leq\, 2|S| \sum_{n= 0}^\infty |D_n| \exp(-a|D_n|) \,\leq\, 2|S|C \sum_{n= 0}^\infty \frac{\log(n+2)}{(n+2)^{Ca}}.
		\end{align*}
		The last expression approaches $0$ as $C \to \infty$, so we can guarantee \eqref{eq:small_sum}.
		
		Consider any $n \in \N$ and $x \in X$. By Lemma~\ref{lemma:concentrated}, we have
		\[
			\P[B_{n,x}] \,\leq\, 2k^{|S|} \exp\left(-\epsilon^2 \frac{|D_n|}{2 |S|^3}\right).
		\]
		If $\dom(B_{n,x}) \cap \dom(B_{m,y}) \neq \0$ for some $m \in \N$ and $y \in X$, then $y \in (SD_m)^{-1}SD_n \cdot x$, and hence for any particular $m \in \N$, there are at most $|S|^2|D_m||D_n|$ choices of such $y$. Therefore, the mapping $\tilde{\omega} \colon \B \to [0;1)$ is a witness to the correctness of $\B$ as long as we have
		\begin{equation}\label{eq:correct1}
			2k^{|S|} \exp\left(-\epsilon^2 \frac{|D_n|}{2 |S|^3}\right) \,\leq\, \omega(n) \prod_{m = 0}^\infty (1 - \omega(m))^{|S|^2 |D_m||D_n|},
		\end{equation}
		for all $n \in \N$. Using the definition of $\omega$ and then taking the logarithm of both sides of \eqref{eq:correct1} and dividing them by $(-|D_n|)$, we rewrite \eqref{eq:correct1} as
		\begin{equation}\label{eq:correct2}
			-\frac{\log(2k^{|S|})}{|D_n|} \,+\, \frac{\epsilon^2}{2 |S|^3} \,\geq\, a \,-\, |S|^2\sum_{m=0}^\infty |D_m|\log(1 - \exp(-a|D_m|)).
		\end{equation}
		Let us first look at the left-hand side of \eqref{eq:correct2}. We have
		\[
			-\frac{\log(2k^{|S|})}{|D_n|} \,+\, \frac{\epsilon^2}{2 |S|^3} \,\geq\, -\frac{\log(2k^{|S|})}{C \log2} \,+\, \frac{\epsilon^2}{2 |S|^3} \,\xrightarrow[C \to \infty]{}\, \frac{\epsilon^2}{2 |S|^3}.
		\]
		As for the right-hand side of \eqref{eq:correct2}, note that $-\log(1 - \xi) < 2\xi$ for all $0 < \xi < 1/2$, so
		\begin{align*}
			a \,-\, |S|^2\sum_{m=0}^\infty &|D_m|\log(1 - \exp(-a|D_m|)) \,<\, a \,+\, 2|S|^2\sum_{m=0}^\infty |D_m| \exp(-a|D_m|) \\
			&\leq\, a \,+\, 2|S|^2 C \sum_{m=0}^\infty \frac{\log(m+2)}{(m+2)^{Ca}} \,\xrightarrow[C \to \infty]{}\, a.
		\end{align*}
		Since $a$ was chosen to be less than $\epsilon^2/(2|S|^3)$, we conclude that \eqref{eq:correct2} holds for all large $C$.
	\end{coolproof}
	
	From \eqref{eq:defect2}, Lemma~\ref{lemma:correct2}, and Theorem~\ref{theo:LLL_meas}, we can derive most of Theorem~\ref{theo:ult_erg_bis}. The only part that is missing is that the map $g \colon \tilde{\Omega}_k \to k$ can be chosen so that 
	\[
		\tilde{\bm{u}}_k(\set{x \in \tilde{\Omega}_k \,:\, g(x) \neq x(\mathbf{1}, 0)}) \,\leq\, \epsilon.
	\]
	To argue this, we have to review the {proof} of Theorem~\ref{theo:LLL_meas}. As mentioned in the introduction, the tool used to prove Theorem~\ref{theo:LLL_meas} is the \emph{Moser--Tardos algorithm}, developed by Moser and Tardos in \cite{MT}. Here we outline only the most relevant elements of the Moser--Tardos theory when applied to our current situation. For further details, see \cite{MT} and \cite[\S3]{MLLL}.
	
	For the rest of this subsection, fix $k \in \N^+$ and a sequence $(\Phi_n)_{n \in \N}$ of bad $k$-events over $\G$. For each $n \in \N$, set $F_n \defeq \dom(\Phi_n)$. Define
	\[
	\B \defeq \B((\Phi_n)_{n \in \N}, \Free(\sigma_{k^\N})), \qquad \text{and} \qquad B_{n,x} \defeq B_x(\Phi_n, \sigma_{k^\N}) \text{ for all } n \in \N \text{ and } x \in \tilde{\Omega}_k.
	\]
	Consider the following inductive construction:
	\begin{leftbar}
		\noindent Set $t_0(x) \defeq 0$ for all $x \in \tilde{\Omega}_k$.
		
		\medskip
		
		\noindent {\sc Step $i \in \N$}: Define
		\begin{align*}
		g_i(x) &\defeq x(\mathbf{1}, t_i(x)) \quad \text{for all } x \in \tilde{\Omega}_k;\\
		\medskip
		A_i' &\defeq \set{(n,x) \in \N \times \tilde{\Omega}_k \,:\, g_i \text{ does not avoid } B_{n,x}}.
		\end{align*}
		Choose $A_i \subseteq A_i'$ to be an arbitrary Borel maximal subset of~$A_i'$ with the property that \[(F_n \cdot x) \cap (F_m \cdot y) = \0 \qquad \text{for all distinct pairs } (n,x),\ (m,y) \in A_i.
		\]
		(Such $A_i$ exists by, e.g., \cite[Lemma 7.3]{KechrisMiller}.) Let
		\[
			T_i \defeq \bigcup_{(n,x) \,\in\, A_i} (F_n \cdot x) \qquad \text{and} \qquad t_{i+1}(x) \defeq \begin{cases}
		t_i(x) + 1 &\text{if } x \in T_i;\\
		t_i(x) &\text{otherwise}.
		\end{cases}
		\]
	\end{leftbar}
	\noindent By definition, $g_0(x) = x(\mathbf{1}, 0)$ for all $x \in \tilde{\Omega}_k$. We call a sequence $\mathcal{A} \defeq (A_i)_{i=0}^\infty$ obtained via the above procedure a \emphd{Borel Moser--Tardos process}. Note that there is not a unique Borel Moser--Tardos process, as there is some freedom in the choice of the Borel maximal subset $A_i \subseteq A_i'$.
	
	Let $\mathcal{A} = (A_i)_{i=0}^\infty$ be a Borel Moser--Tardos process. For $x \in \tilde{\Omega}_k$, define $t(x) \in \N \cup \set{\infty}$ by
	\[
		t(x) \defeq \lim_{i \to \infty} t_i(x).
	\]
	We say that $x$ is \emphd{$\mathcal{A}$-stable} if $t(x) < \infty$, i.e., if the corresponding sequence $t_0(x)$, $t_1(x)$, \ldots{} is eventually constant. Let $\Stab(\mathcal{A}) \subseteq \tilde{\Omega}_k$ denote the set of all $\mathcal{A}$-stable elements. For $x \in \Stab(\mathcal{A})$, we can define
	\begin{equation}\label{eq:g}
		g(x) \defeq x(\mathbf{1}, t(x)).
	\end{equation}
	It is easy to verify \ep{see \cite[Proposition 3.3]{MLLL}} that if $F_n \cdot x \subseteq \Stab(\mathcal{A})$, then $x \not \in \Def(g, \Phi_n, \sigma_{k^\N})$. The \emphd{index} $\Ind(n, x, \mathcal{A}) \in \N \cup \set{\infty}$ of a pair $(n,x) \in \N \times \tilde{\Omega}_k$ in $\mathcal{A}$ is defined by the formula
	\[
	\Ind(n, x, \mathcal{A}) \defeq |\set{i \in \N\,:\, (n,x) \in A_i}|.
	\]
	Note that for all $x \in \Free(\tilde{\Omega}_k)$, we have
	\begin{equation}\label{eq:stabilizing1}
	t(x) \,=\, \sum_{n = 0}^\infty \,\sum_{\delta \in F_n} \Ind(n, \delta^{-1} \cdot x, \mathcal{A}),
	\end{equation}
	and hence such $x$ is $\mathcal{A}$-stable if and only if the expression on the right hand side of \eqref{eq:stabilizing1} is finite. The following theorem is the central result of the Moser--Tardos theory:
	
	\begin{theo}[{Moser--Tardos~\cite{MT}; see also \cite[Theorem~3.5]{MLLL}}]\label{theo:MoserTardos1}
		Let $\omega \colon \N \to [0;1)$ be a function such that the mapping $\tilde{\omega} \colon \B \to [0;1) \colon B_{n,x} \mapsto \omega(n)$ is a witness to the correctness of $\B$. Then, for any Borel Moser--Tardos process $\mathcal{A}$ and for all $n \in \N$, we have
		\begin{equation*}
		\int_{\tilde{\Omega}_k} \Ind(n, x, \mathcal{A}) \, \Diff \tilde{\bm{u}}_k(x) \,\leq\, \frac{\omega(n)}{1 - \omega(n)}.
		\end{equation*}
	\end{theo}

	\begin{corl}[to Theorem~\ref{theo:MoserTardos1}]\label{corl:dist}
		Let $\omega \colon \N \to [0;1)$ be such that  $\tilde{\omega} \colon \B \to [0;1) \colon B_{n,x} \mapsto \omega(n)$ is a witness to the correctness of $\B$. Then there is a Borel function $g \colon \tilde{\Omega}_k \to k$ such that
		\[
			\tilde{\bm{u}}_k(\Def(g, (\Phi_n)_{n \in \N}, \sigma_{k^\N})) = 0
		\qquad
		\text{and}
		\qquad
			\tilde{\bm{u}}_k(\set{x \in \tilde{\Omega}_k \,:\, g(x) \neq x(\mathbf{1}, 0)}) \,\leq\, \sum_{n=0}^\infty |F_n| \cdot \frac{\omega(n)}{1 - \omega(n)}.
		\]
	\end{corl}
	\begin{coolproof}
		First we show that the sum
		\[
		S \,\defeq\, \sum_{n=0}^\infty |F_n| \cdot \frac{\omega(n)}{1 - \omega(n)}
		\]
		is finite. Without loss of generality, assume that $\Phi_0 \neq \0$. Consider any $x \in \Free(\tilde{\Omega}_k)$. Since $\tilde{\omega}$ is a witness to the correctness of $\B$, we have $\P[B_{0, x}] \leq \omega(0) < 1$, so $F_0 \neq \0$. Hence, for every $n \in \N^+$, there exist at least $|F_n|$ distinct $y$ with $B_{n,y} \in \Nbhd_\B(B_{0,x})$. Therefore,
		$
			\prod_{n=1}^\infty (1 - \omega(n))^{|F_n|} \geq \P[B_{0,x}] > 0,
		$
		which implies that $\sum_{n=0}^\infty |F_n| \omega(n)$ is finite. In particular, for all sufficiently large $n$ we have $\omega(n) \leq 1/2$ and $\omega(n)/(1- \omega(n)) \leq 2\omega(n)$, and hence $S$ is also finite.
		
		Let $\mathcal{A} = (A_i)_{i=0}^\infty$ be an arbitrary Borel Moser--Tardos process and let $g$ be given by \eqref{eq:g}. From \eqref{eq:stabilizing1} and the \hyperref[theo:MoserTardos1]{Moser--Tardos theorem}, we get
		\begin{align*}
		\int_{\tilde{\Omega}_k} t(x) \,\Diff \tilde{\bm{u}}_k(x) \,&=\, \sum_{n = 0}^\infty \sum_{\delta \in F_n} \int_{\tilde{\Omega}_k} \Ind(n, \delta^{-1} \cdot x, \mathcal{A}) \, \Diff \tilde{\bm{u}}_k(x) \\
		[\text{$\tilde{\bm{u}}_k$ is shift-invariant}] \qquad \,&=\, \sum_{n = 0}^\infty |F| \cdot \int_{\tilde{\Omega}_k} \Ind(n, x, \mathcal{A}) \, \Diff \tilde{\bm{u}}_k(x) \,\leq\, S\,<\,\infty.
		\end{align*}
		In particular, $t(x) < \infty$ for $\tilde{\bm{u}}_k$-a.e.\ $x \in \tilde{\Omega}_k$, i.e., $\tilde{\bm{u}}_k(\Stab(\mathcal{A})) = 1$, so
		\[
			\tilde{\bm{u}}_k(\Def(g, (\Phi_n)_{n \in \N}, \sigma_{k^\N})) = 0.
		\]
		Furthermore, if $x \in \Stab(\mathcal{A})$ and $g(x) \neq x(\mathbf{1}, 0) = g_0(x)$, then $t(x) \geq 1$; thus,
		\[
		\tilde{\bm{u}}_k(\set{x \in \tilde{\Omega}_k \,:\, g(x) \neq x(\mathbf{1}, 0)}) \,\leq\, \tilde{\bm{u}}_k(\set{x \in \tilde{\Omega}_k \,:\, t(x) \geq 1}) \,\leq\, \int_{\tilde{\Omega}_k} t(x) \,\Diff \tilde{\bm{u}}_k(x) \,\leq\,  S,
		\]
		as desired.
	\end{coolproof}
	
	Since the domain of $\Phi(k, S, \epsilon, D)$ is, by definition, $SD$, \eqref{eq:small_sum} in the statement of Lemma~\ref{lemma:correct2} and Corollary~\ref{corl:dist} yield the remaining part of Theorem~\ref{theo:ult_erg_bis}.

	\section{Open problems}\label{sec:remks}
	
	The following is perhaps the central open question regarding the behavior of the LLL in the measurable setting:
	
	\begin{prob}\label{prob:meas_LLL}
		Does the SLLL hold measurably with a null defect? In other words, can one replace $\mu(\Def(f, \B)) \leq \delta$ by $\mu(\Def(f, \B)) = 0$ in the conclusion of Theorem~\ref{theo:LLL_approx}?
	\end{prob}
	
	\noindent A positive solution to Problem~\ref{prob:meas_LLL} would allow one to strengthen Theorem~\ref{theo:main_delta} by taking $\delta=0$. For now, we leave this potential strengthening as an open problem.
	
	\begin{prob}\label{prob:delta_zero}
		Does Theorem~\ref{theo:main_delta} hold with $\delta = 0$?	
	\end{prob}
	
	\noindent As mentioned in \S\ref{subsec:MLLL}, the SLLL {fails} in the purely Borel context~\cite[Theorem~1.6]{CJMST-D}. However, it is still conceivable that a purely Borel pointwise version of the \hyperref[theo:AW]{Ab\'ert--Weiss theorem}, similar to Theorem~\ref{theo:main_Borel}, holds in full generality, in which case a different proof approach might be needed to establish it. We state it here as another open question.
	
	\begin{prob}\label{prob:Borel}
		Let $K$ be a compact metric space and let $f \colon \Omega \to K$ be a Borel function. Fix an open neighborhood $U$ of the measure $\M_\LG \code{f}$. Does there always exist a nonempty finite set $D \subset \G$ such that the following statement holds?
		
		\smallskip
		
		Let $\alpha \colon \G \acts X$ be a free Borel action of $\G$. Then there is a Borel map $g \colon X \to K$ such that \[\M_D \code{g} (x) \in U, \qquad \text{for all}\ x \in X.\]
	\end{prob}

		\printbibliography

	\appendix
	
	\section{Proof of Theorem~\ref{theo:bad}}\label{sec:app}
	
	Let $\alpha \colon \Z \acts (X, \mu)$ be a free \pmp action of $\Z$. For a Borel set $A \subseteq X$, let $[A]$ denote the class of all Borel sets $B \subseteq X$ with $\mu(A \symdif B) = 0$. The measure algebra $\mathrm{MAlg}(X, \mu)$ is the space of all classes $[A]$ with the metric $\dist([A], [B]) \defeq \mu(A \symdif B)$. Note that the space $\mathrm{MAlg}(X, \mu)$ is Polish. For a sequence $(D_n)_{n \in \N}$ of nonempty finite subsets of $\G$, let
	\[
		\mathfrak{L}(\alpha, X, \mu, (D_n)_{n \in \N}) \defeq \set{[A] \in \mathrm{MAlg}(X, \mu) \,:\, \textstyle \liminf_{n \to \infty} \E_{D_n} \mathbbm{1}_A(x) = 0 \text{ for $\mu$-a.e.}\ x \in X};
	\]
	\[
	\mathfrak{U}(\alpha, X, \mu, (D_n)_{n \in \N}) \defeq \set{[A] \in \mathrm{MAlg}(X, \mu) \,:\, \textstyle \limsup_{n \to \infty} \E_{D_n} \mathbbm{1}_A(x) = 1 \text{ for $\mu$-a.e.}\ x \in X}.
	\]
	It is straightforward to check that the sets $\mathfrak{L}(\alpha, X, \mu, (D_n)_{n \in \N})$ and $\mathfrak{U}(\alpha, X, \mu, (D_n)_{n \in \N})$ are $G_\delta$ in $\textrm{MAlg}(X, \mu)$. Therefore, to establish the conclusion of Theorem~\ref{theo:bad}, it is enough to ensure that both these sets are dense. Below we only give the argument that shows that $\mathfrak{U}(\alpha, X, \mu, (D_n)_{n \in \N})$ is dense; the proof for $\mathfrak{L}(\alpha, X, \mu, (D_n)_{n \in \N})$ is the same, \emph{mutatis mutandis}.

	\begin{lemma}\label{lemma:epsilon}
		Let $h \colon \N \to \N$ be an arbitrary function and let $\epsilon > 0$. There exists a finite sequence $D_0$, \ldots, $D_{N-1}$ of finite subsets of $\Z$ with the following properties:
		\begin{itemize}[label=--,wide]
			\item each $D_n$ is an interval;
			
			\item $|D_n| \geq h(n)$ for all $0 \leq n < N$;
			
			\item for every free \pmp action $\Z \acts (X, \mu)$, there is a Borel set $A \subseteq X$ with $\mu(A) \leq \epsilon$ such that
			\[
			\mu(\set{x \in X \,:\, D_n \cdot x \subseteq A \text{ for some } 0 \leq n < N}) \,\geq 1-\epsilon.
			\]
		\end{itemize}
	\end{lemma}
	\begin{coolproof}
		Take any $N$ so large that
		\[
			\frac{2}{N+1} \,<\, \epsilon \qquad \text{and} \qquad (1-\epsilon/2)\frac{N}{N+1} \,>\, 1- \epsilon.
		\]
		Let $\ell \defeq \max_{n=0}^{N-1} h(n)$ and for each $0 \leq n < N$, define
		\[
			D_n \defeq \set{n\ell, \, n\ell+1, \, \ldots, \, n\ell + \ell - 1}.
		\]
		We claim that this sequence of intervals works. Let $\alpha \colon \Z\acts (X, \mu)$ be a free \pmp action of $\Z$ induced by a measure\-/preserving transformation $T \colon X \to X$. By Rokhlin's lemma, there exists a Borel set $R \subseteq X$ such that its translates $R$, $TR$, \ldots, $T^{(N+1)\ell-1} R$ are pairwise disjoint and their union has measure at least $1 - \epsilon/2$. Let
		\[
			A \defeq \bigcup_{i = (N-1)\ell}^{(N+1)\ell-1} T^i R \qquad \text{and} \qquad B \defeq \bigcup_{i = 0}^{N\ell-1} T^{i} R.
		\]
		\begin{figure}[H]
			\begin{tikzpicture}
				\draw (0,0) -- (10,0);
				\node[anchor=east] at (0,0) {\LARGE\ldots};
				\node[anchor=west] at (10,0) {\LARGE\ldots};
				
				\draw (1,-0.2) -- (1,0.2);
				\draw (2,-0.2) -- (2,0.2);
				\draw (3,-0.2) -- (3,0.2);
				\draw (4,-0.2) -- (4,0.2);
				\draw (5,-0.2) -- (5,0.2);
				\draw (6,-0.2) -- (6,0.2);
				\draw (7,-0.2) -- (7,0.2);
				\draw (8,-0.2) -- (8,0.2);
				\draw (9,-0.2) -- (9,0.2);
				
				\draw [decorate,decoration={brace,mirror,amplitude=10pt}] (1,-0.25) -- (8,-0.25) node [midway,anchor=north,yshift=-10pt] {$B$};
				
				\draw [decorate,decoration={brace,amplitude=10pt}] (7,0.25) -- (9,0.25) node [midway,anchor=south,yshift=10pt] {$A$};
			\end{tikzpicture}
			\caption{A cartoon of the sets $A$ and $B$.}\label{fig:Rokhlin}
		\end{figure}
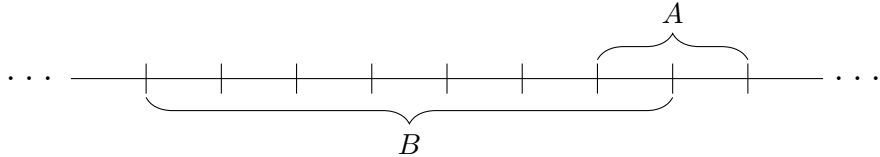
		\noindent Then $\mu(A) = 2\ell \mu(R) \leq 2/(N+1) < \epsilon$ and $\mu(B) = N\ell\mu(R) \geq (1-\epsilon/2) N/(N+1) > 1 - \epsilon$, and for each $x \in B$, there is some $0 \leq n < N$ with $D_n \cdot x \subseteq A$, as desired (see Fig.~\ref{fig:Rokhlin}).
	\end{coolproof}
	
	Let $h \colon \N \to \N$ be any function. Applying Lemma~\ref{lemma:epsilon} repeatedly, we construct an increasing sequence of natural numbers $(N_i)_{i \in \N}$ starting with $N_0 \defeq 0$ and a sequence of finite intervals $(D_n)_{n \in \N}$ with $|D_n| \geq h(n)$, such that for every free \pmp action $\Z \acts (X, \mu)$, there exists a sequence of Borel sets $(A_i)_{i \in \N}$ with $\mu(A_i) \leq 2^{-i-1}$ and $\mu(B_i) \,\geq 1-2^{-i-1}$, where
	\[
		B_i \defeq \set{x \in X \,:\, D_n \cdot x \subseteq A_i \text{ for some } N_i \leq n < N_{i+1}}.
	\]
	For $k \in \N$, let $A_{\geq k} \defeq \bigcup_{i = k}^\infty A_i$. We claim that $A_{\geq k} \in \mathfrak{U}(\alpha, X, \mu, (D_n)_{n \in \N})$. Indeed, \[\limsup_{n \to \infty} \E_{D_n} \mathbbm{1}_{A_{\geq k}}(x) = 1 \qquad \text{for all } x \in \limsup_{i \to \infty} B_i,\] and, by Fatou's lemma, $\mu(\limsup_{i \to \infty} B_i) \geq \limsup_{i \to \infty} \mu(B_i) =1$. Now if $[S] \in \mathrm{MAlg}(X, \mu)$, then $[S \cup A_{\geq k}] \in \mathfrak{U}(\alpha, X, \mu, (D_n)_{n \in \N})$ as well, and $\dist([S], [S \cup A_{\geq k}]) \leq \mu(A_{\geq k}) \leq \sum_{i=k}^\infty \mu(A_i) \leq 2^{-k}$. Since $k$ is arbitrary, this shows that $\mathfrak{U}(\alpha, X, \mu, (D_n)_{n \in \N})$ is dense in $\mathrm{MAlg}(X, \mu)$, as desired.

\end{document}